# A study of the Structural Properties of finite $G$-graphs and their Characterisation.

By

Lord Clifford Kavi

(10280594)

lord@aims.edu.gh
June 2016

THIS THESIS IS SUBMITTED TO THE UNIVERSITY OF GHANA, LEGON IN PARTIAL FULFILMENT OF THE REQUIREMENTS FOR THE AWARD OF MASTER OF PHILOSOPHY MATHEMATICS DEGREE.

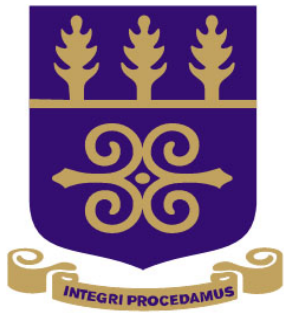

# DECLARATION

This thesis was written in the Department of Mathematics, University of Ghana, Legon from August 2015 to July 2016 in partial fulfilment of the requirements for the award of Master of Philosophy degree in Mathematics under the supervision of Dr. Margaret McIntyre and Dr. Douglass Adu-Gyamfi of the University of Ghana.

I hereby declare that except where due acknowledgement is made, this work has never been presented wholly or in part for the award of a degree at the University of Ghana or any other University.

Signature: 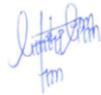
Student: Lord Clifford Kavi
lord@aims.edu.gh

Signature: .................................................
Dr. Margaret McIntyre

Signature: .................................................
Dr. Douglass Adu-Gyamfi

# ABSTRACT

The $G$-graph $\Gamma(G, S)$ is a graph from the group $G$ generated by $S \subseteq G$, where the vertices are the right cosets of the cyclic subgroups $\langle s \rangle, s \in S$ with $k$-edges between two distinct cosets if there is an intersection of $k$ elements. In this thesis, after presenting some important properties of $G$-graphs, we show how the $G$-graph depends on the generating set of the group. We give the $G$-graphs of the symmetric group, alternating group and the semi-dihedral group with respect to various generating sets. We give a characterisation of finite $G$-graphs; in the general case and a bipartite case. Using these characterisations, we give several classes of graphs that are $G$-graphs. For instance, we consider the Turán graphs, the platonic graphs and biregular graphs such as the Levi graphs of geometric configurations. We emphasis the structural properties of $G$-graphs and their relations to the group $G$ and the generating set $S$.

As preliminary results for further studies, we give the adjacency matrix and spectrum of various finite $G$-graphs. As an application, we compute the energy of these graphs. We also present some preliminary results on infinite $G$-graphs where we consider the $G$-graphs of the infinite group $SL_2(\mathbb{Z})$ and an infinite non-Abelian matrix group.

**DEDICATION**

To my lovely daughter, Charis.

# ACKNOWLEDGEMENTS

"Surely, God is my help; the Lord is the one who sustains me. Psalm 54:4".

I foremost want to thank my Lord and Saviour Jesus Christ for His grace that sustains me.

I want to appreciate African Institute for Mathematical Sciences (AIMS-Ghana), for supporting this research financially.

I also want to express my sincere gratitude to my main supervisor, Dr. Margaret McIntyre. I really appreciate your time, resources and kindness to me. I have enjoyed working with you. Also to my other supervisors Dr. Douglass Adu-Gyamfi and Prof. Nancy Ann Neudauer, who gave me relevant suggestions.

My sincere gratitude also goes to all the teaching and non teaching staff and colleague students of the department. You made our stay here worthwhile.

To them who though never had classmates yet gave me the opportunity to have countless class-mates. I love you Mum and Dad, Mr. and Mrs. Vicent Kavi.

Though I am your son-in-law, you love me and treat me as a biological son, Rev. Prof. and Mrs. James Ansah, God bless you.

To my best friend who doubles as my beloved wife, Rhoda, thanks for your patience and understanding. You define love.

# Contents

# Chapter 1

# Introduction and Motivation

Graphs constructed from groups are useful in applications such as optimization on parallel architectures, or for the study of interconnection networks because of their regularity and their underlying algebraic structure [11]. One such graph is the Cayley graph. Given a group $G$ generated by $S = \{s_1, s_2, ..., s_k\}$, the Cayley graph $(G, S)$ is the graph with vertex set $V$ consisting of the elements of $G$. An edge connects vertices $g_1, g_2 \in V$ if there exists some $s_i \in S$ such that $g_1 = s_i g_2$. Cayley graphs have been extensively studied, see [4, 18]. Some interesting properties of Cayley graphs are that they are symmetric, regular, and vertex-transitive, properties which make them very useful, but also present some limitations in their application. Cayley graphs are useful in solving rearrangement problems and the design of interconnection networks for parallel CPU's [18].

In 2005, Bretto and Faisant [7] introduced and defined a new type of graph, which they called a *G-graph*, a graph from group. Given a group $G$ generated by $S = \{s_1, s_2, ..., s_k\}$, the $G$-graph $\Phi(G, S)$ was defined as the graph whose vertices are the cycles $(x, sx, ..., s^{o(s)-1}x), x \in G$ of the left $s$-translation on $G$. The $G$-graph has $p$-edges between any two cycles with an intersection of $p$ elements. The $G$-graph defined this way has $o(s)$ loops on the vertex $(x, sx, ..., s^{o(s)-1}x)$. Later, in [8, 11, 12], they defined the $G$-graph without loops. In 2013, [35], Tanasescu, et al gave a variation of the definition for $G$-graphs, where the vertices are the right cosets of the cyclic subgroups $\langle s \rangle$. There exist $k$-edges between two distinct cosets which have an intersection of $k$ elements.

Like Cayley graphs, $G$-graphs have specific symmetries and regularity. So they can be used in any area of science where Cayley graphs can be used. Because unlike Cayley graphs, $G$-graphs can be semi-regular, edge-transitive and $G$-graphs are $|S|$-partite. Thus they are very useful in constructing symmetric and semi-symmetric networks [13]. Various characteristics of the group can be seen on the corresponding $G$-graph [11]. See [34] for some of the properties of the group that can be observed on the graph. These properties of a $G$-graph mean that they can be applied to problems where the limited information from a Cayley graph is insufficient to solve the problem. Just like Cayley graphs, $G$-graphs are very regular and have compact representations from groups, and these motivate their study for the analysis and design of networks and applications in areas such as communication networks, code theory and cryptography [34]. It is worth knowing that, there exists a relation between Cayley graphs and $G$-graphs. In fact in [9], the authors gave a sufficient condition to recognize when a $G$-graph is a Cayley graph. It is also shown that the line graph of a bipartite $G$-graph is a Cayley graph [13].

In [13], the authors gave a characterisation of bipartite $G$-graphs. The authors of [22] also gave a characterisation of $G$-graphs with and without loops in terms of the automorphism on the graph. It is shown in [19] that any connected bipartite $G$-graph, and the $G$-graph of a group with symmetric presentation, are optimally connected, a result that shows how networks constructed from such $G$-graphs can be very robust. It is known that Cayley graphs are optimally connected because they are vertex-transitive [21]. Semi-symmetric $G$-graphs provide a better option than Cayley graphs when the network needs not be symmetric and totally regular.

This thesis seeks to understand the structure of $G$-graphs, and to contribute to the recognition and characterization problem. We also consider infinite $G$-graphs. This thesis is organised in chapters as follows:

1. Chapter 1, the current chapter provides the definition and the motivation for the study of $G$-graphs.

2. Chapter 2 deals with the necessary definitions we need for understanding the rest of the document. We also make some useful contributions toward the properties of $G$-graphs.

3. Chapter 3 shows how the $G$-graph depends on the generating set of the group.

4. Chapter 4 presents results on the structure of the $G$-graph of a symmetric group and that of an alternating group generated by various generating sets.

5. Chapter 5 highlights the $G$-graphs of the quaternions, the dihedral and the semi-dihedral groups.

6. In Chapter 6, we present results that deal with the $G$-graph recognition problem. We give necessary and sufficient conditions for a given graph to be the $G$-graph $\Gamma(G, S)$ of some group $G$ with generating set $S$. As application, we describe graphs which are $G$-graphs and others which are not. For example, we show which Turán graphs are $G$-graphs and which are not, we look at the platonic graphs, the $n$-dimensional hypercube, and we give a characterisation of bipartite $G$-graphs.

7. In chapter 7, we present some preliminary results on the adjacency matrix and the spectrum of a finite $G$-graph. As an application, we consider the energy of some finite $G$-graphs. We then shift our focus to infinite $G$-graphs. We show the $G$-graph of some infinite groups that are finitely generated. These are important aspects of $G$-graphs which will require further investigation.

8. The last chapter deals with our conclusion and presents some open questions. We also give some equally important aspect of a $G$-graph that remains for future work.

## Chapter 2

## Definitions and Properties

This chapter looks at some concepts in group theory and graph theory.

## 2.1 Group Theory Definitions

For the prerequisites from group theory we refer to [24, 30] and [31].

**Definition 2.1.1.** A *group* $G$ is a set with a binary operation $*$, satisfying the following axioms:

- For all $a, b, c \in G$, $a * (b * c) = (a * b) * c$. That is, the group operation is associative.

- There exists an identity $e \in G$ such that for all $a \in G$, $a * e = e * a$.

- For all $a \in G$, there exists $a^{-1}$ called an inverse of $a$ such that
  $a * a^{-1} = a^{-1} * a = e$.

**Definition 2.1.2.** The *order* $|G|$ of the group $G$ is its cardinality. If $G$ has a finite order, we say $G$ is finite otherwise it is infinite. The order $o(a)$ of $a \in G$ is the minimum $k \in \mathbb{N}$ such that $a^k = e$. If there exists no such $k$ then, we say $a$ has infinite order.

**Definition 2.1.3.** *Generating set*: Let $G$ be a group and $S$ a subset of $G$. $S$ is a *generating set* of $G$ written as $G = \langle S \rangle$ if every element of $G$ can be written as a product of elements of $S$ and inverses of elements of $S$.

*Cyclic group*: A group $G$ is a *cyclic group* if it can be generated by a single element say $g$. Then $G = \langle g \rangle$.

**Definition 2.1.4.** *Subgroup*: A subset $H$ of a group $G$, which inherits the group structure of $G$ is a *subgroup* of $G$. That is:

- for all $a, b \in H$, $ab \in H$, so closure is satisfied.
- the identity $e \in G$ is contained in $H$.
- for all $a \in H$, there exists an inverse $a^{-1} \in H$.

*Cyclic subgroup*: A subgroup that can be generated by a single element is called a *cyclic subgroup*. For instance, the cyclic subgroup of $a \in G$ is $\{e, a, a^2, a^3, ..., a^{o(a)-1}\}$.

**Definition 2.1.5.** *Homomorphism*: A *homomorphism* from a group $G$ to a group $G'$ is a function $f : G \to G'$ such that $f(ab) = f(a)f(a)$ for all $a, b \in G$.

Two groups $G$ and $H$ are *isomorphic* written as $G \cong H$ if there exists a bijective homomorphism $f : G \to H$. The map $f$ is an *isomorphism*.

Suppose $S$ generates $G$ and $S'$ generates $G'$. Then $(G, S) \cong (G', S)$ if there exists an isomorphism $f : G \to G'$ such that $f(S) = S'$.

**Definition 2.1.6.** *Group presentation*: A group presentation is a definition or construction of a group using generators and relations between the generators. $\langle S | R \rangle$ is a presentation of the group generated by the generating set $S$ subject to every relation $(u, v) \in R$.

**Definition 2.1.7.** *Right coset*: Let $H$ be a subgroup of $G$. For $x \in G$, the set $Hx = \{hx | h \in H\} \subset G$ is called the *right coset* of $H$ containing $x$. The right cosets of $H$ form a partition of $G$. That is, for any $x, y \in G$, either $Hx = Hy$ or $Hx \cap Hy = \emptyset$. The left coset is defined similarly. In this thesis, anytime we use coset, we mean right coset.

## 2.2 Graph Theory Definitions

For the prerequisites from graph theory, we refer to [6] and [26].

**Definition 2.2.1.** A *graph* $\Gamma$ is a triple $(V, E, \tau)$ where $E$ is an edge set, $V$ is a vertex set and $\tau$ is a map that associates to each edge $e \in E$ an unordered pair of vertices $a, b \in V$. If $\tau(e) = \{a, b\}$, then edge $e$ is said to be incident to $a$ and $b$. If $a$ and $b$ are not distinct, we say $e$ is a *loop*. If more than one edge is incident to the same pair of vertices, then we have *multiple edges*. A graph that has neither multiple edges nor loops is a *simple graph*. We shall strictly consider graphs without loops. Thus all our results apply to graphs without loops. A graph is *finite* if $|E| \neq \infty \neq |V|$ otherwise it is *infinite*. Unless stated otherwise, we deal with *finite* graphs. Graphs that have orientation are called *directed graphs*. In this thesis, all our graphs are *undirected*.

**Definition 2.2.2.** The *degree* of a vertex $v$ is the number of edges that are incident to it, with loops counting twice. We shall denote the degree of a vertex $v$ as $deg(v)$. If every vertex of a graph has the same degree say $d$, we say the graph is *d-regular*.

**Definition 2.2.3.** *k-partite graph*: A graph $\Gamma = (V, E, \tau)$ is said to be *k*-partite if we can partition the vertex set $V$ into $V = V_1 \cup \cdots \cup V_k$, $V_i \cap V_j = \emptyset$ such that there is no edge between vertices in the same partition. In other words, the graph is *k*-chromatic, the minimum number of colours needed to colour the vertices of a graph in such a way that, no adjacent vertices have the same colour is $k$. When $k = 2$, we say the graph is bipartite, that is $\Gamma = (V_1 \cup V_2, E)$.

*Semi-regular graph*: A graph is said to be *semi-regular* if vertices in the graph assume one of two degrees. Thus a $(x, y)$-semi-regular graph is such that, for all $v_I \in V$, $deg(v_i) = x$ and $deg(v_j) = y$.

*Regular-in-partition*: We say a $k$-partite graph is *regular-in-partition* if vertices in the same partition have the same degree. For instance, a semi-regular bipartite graph also known as a biregular graph has degree of vertices in the same bipartition to be the same. Thus, biregular graphs are regular-in-partition.

**Definition 2.2.4.** *Subgraph*: Let $G = (V, E, \tau)$ be a graph.

The graph $H = (V', E', \tau)$ is a *subgraph* of $G$ if $V' \subset V$ and $E' \subset E$.

*Induced subgraph*: If $V' \subset V$. The subgraph of $G$ induced by the vertex set $V'$ is the graph $H = (V', E', \tau)$ where $\{v_1, v_2\} \in E$ if and only if $v_1$ and $v_2$ are both in $V'$.

**Definition 2.2.5.** *Graph isomorphism*: Two graphs $H$ and $H^*$ are isomorphic if there exists a bijection $f : V(H) \to V(H^*)$ such that $f(u)$ and $f(v)$ are adjacent in $H^*$ if and only if $u$ and $v$ are adjacent in $H$. That is, the bijection between their vertex sets preserves adjacency. An isomorphism of a graph onto itself is called an *automorphism*. The set of automorphisms of a graph $\Gamma$ form a group. We shall denote the automorphism group as $Aut(\Gamma)$.

## 2.3 $G$-graph Definition

The idea of a $G$-graph was introduced and defined in [7]. A slight variation to the definition was given in [35]. We follow the definition given in [35].

**Definition 2.3.1.** Let $G$ be a group and $S$ a set of elements of $G$ such that $G = \langle S \rangle$. The $G$-graph $\Phi(G, S)$ is the graph $\Gamma = (V, E, \tau)$ such that:

1. The set of vertices of $\Phi(G, S)$ is $V = \bigcup_{s \in S} V_s$, where $V_s = \{\langle s \rangle x, x \in G\}$ and $\langle s \rangle x$ is the right coset of the cyclic subgroup $\langle s \rangle$ containing $x$.

2. For $\langle s \rangle x, \langle t \rangle y \in V, (s \neq t)$, there exist $k$ edges between $\langle s \rangle x$ and $\langle t \rangle y$ if $|\langle s \rangle x \cap \langle t \rangle y| = k$.

By this definition, $\Phi(G, S)$ is a graph with $o(s)$ loops at vertex $\langle s \rangle x$.

Throughout this thesis, we shall deal with $G$-graphs without loops. We shall denote the $G$-graph without loops by $\Gamma(G, S)$.

## 2.4 Properties of finite $G$-graphs

In this section, we deal with finite groups and finite $G$-graphs. We shall denote a group by $G$ and a generating set by $S$.

**Proposition 2.4.1.** *[10] If $|S| = n$, then the $G$-graph $\Gamma(G, S)$ is $n$-partite .*

The following proposition also appears in [34].

**Proposition 2.4.2.** *The degree of any vertex $v$ in the $G$-graph $\Gamma(G, S)$ of a group $G$ is $deg(v) = o(s)(|S| - 1)$ where $o(s)$ is the order of the generator $s$ and $|S|$ is the cardinality of the generating set $S$.*

*Proof.* For any group $G$ and for any generating set $S$ of $G$, $\Gamma(G, S)$ is an $|S|$-partite graph (see proposition 2.4.1.)

Suppose $|S| = n$, then we have the vertex set $V$ of $\Gamma(G, S)$ as $V = V_1 \cup \ldots \cup V_n$ where the vertices in the partition $V_i$ are the cosets $\langle s_i \rangle g$, for $s_i \in S$ and $g \in G$ with $G = \bigcup_{g \in G} \langle s_i \rangle g$ (in fact the vertices in $V_i$ partition the elements of $G$). Since vertices $\langle s_i \rangle\, g$, $\langle s_j \rangle\, h$ have an edge between them if and only if $\langle s_i \rangle\, g \cap \langle s_j \rangle\, h \neq \emptyset$, this implies $\Gamma$ is partite and each vertex has degree $deg(v) = o(s_i)(|S|-1)$. The vertex (as a set) contains $o(s_i)$ elements of $G$, each of which appear in exactly one vertex of each other set $V_j$, $j \neq i$, of the partition of vertices. □

The following are immediate from proposition 2.4.2.

**Corollary 2.4.3.** *There are $\frac{|G|}{o(s_i)}$ distinct cosets labelling the vertices in partition $V_i$ corresponding to cosets of $G$ by $\langle s_i \rangle$. So if $|S| = n$, there are $\sum_{i=1}^{n} \frac{|G|}{o(s_i)}$ vertices in $\Gamma(G,S)$.*

**Corollary 2.4.4.** *If $|S| = n$, then the number of edges in $\Gamma(G,S)$ is $\frac{n(n-1)}{2}|G|$.*

*Proof.* The number of edges will be $\frac{1}{2}\sum_{i=1}^{n} \frac{|G|}{o(s_i)} \times o(s)(|S|-1)$. If $|S| = n$, then $\frac{1}{2}\sum_{i=1}^{n} \frac{|G|}{o(s_i)} \times o(s)(|S|-1) = \frac{1}{2}n|G|(n-1) = \frac{n(n-1)}{2}|G|$. □

**Proposition 2.4.5.** *The G-graph $\Gamma(G,S)$ of the group $G$ with the generating set $S$ is regular in each vertex set of the partition $V = V_1 \cup ... \cup V_{|S|}$.*

*Proof.* Let the generating set $S$ of $G$ be $S = \{s_1, s_2, ..., s_n\}$ and let the order of each $s_i \in S$ be $o(s_i) = k_i$, not necessarily distinct. Since $|S| = n$, we have the vertex set $V$ of $\Gamma(G,S)$ as $V = V_1 \cup ... \cup V_n$, where vertices in the partition $V_i$ are the cosets $\langle s_i \rangle\, g$, for $s_i \in S$ and $g \in G$ with $G = \bigcup_{g \in G} \langle s_i \rangle\, g$. The cardinality of a coset $\langle s_i \rangle\, g$ is $o(s_i) = k_i$. Thus all cosets of the cyclic subgroup $\langle s_i \rangle$ have the same cardinality, $k_i$. Then by Prop. 2.4.2, the degree of each vertex in a vertex set $V_i$ is $deg(v) = k_i(n-1)$. Thus vertices in the same vertex set of the partition of vertices have the same degree. □

**Corollary 2.4.6.** *If $|S| = n$ and $o(s_i) = k$ for all $s_i \in S$, then the G-graph $\Gamma(G,S)$ of the group $G$ with the generating set $S$ is regular of degree $k(n-1)$.*

*Proof.* From Prop. 2.4.5, the degree of each vertex in a vertex set (partition set) is dependent on the order of the corresponding generator $s_i$. Since $o(s_i) = k$ for all $s_i \in S$, each vertex set has degree $k(n-1)$. Thus $\Gamma(G,S)$ is regular. □

The following are also immediate from proposition 2.4.2.

**Corollary 2.4.7.** *If $|S|$ is odd for all $|S| > 2$, then the G-graph $\Gamma(G,S)$ is Eulerian.*

*Proof.* Any vertex $v$ has degree $deg(v) = o(s)(|S| - 1)$, so if $|S|$ is odd, then $(|S| - 1)$ is even and hence $deg(v)$ is even for all $v$. □

**Corollary 2.4.8.** *If $|S| = 2$, then the G-graph $\Gamma(G, S)$ is Eulerian if and only if all of its generators are of even order, (that is, $o(s_1) = 2x$ and $o(s_2) = 2y$, $x, y \in \mathbb{N}$).*

*Proof.* Suppose $|S| = 2$ and let $S = \{s_1, s_2\}$. Then $\Gamma(G, S)$ is bipartite by prop 2.4.1. Suppose $\Gamma(G, S)$ is Eulerian. $\Gamma(G, S)$ is regular in each vertex set of the partition (see Prop. 2.4.5), so vertices in one vertex set have the degree $deg = o(s_1)(2 - 1) = o(s_1)$ while the vertices in the other vertex set have the degree $deg = o(s_2)(2 - 1) = o(s_2)$. Since it is Eulerian, each vertex has even degree so $o(s_1) = 2x$ and $o(s_2) = 2y$, $x, y \in \mathbb{N}$.

Conversely, suppose $o(s_1) = 2x$ and $o(s_2) = 2y$, $x, y \in \mathbb{N}$. This means the cosets of the cyclic subgroup $\langle s_1 \rangle$ have cardinality $2x$ and the cosets of the cyclic subgroup $\langle s_2 \rangle$ have cardinality $2y$. Thus, the degree of vertices in the vertex set corresponding to cosets of the cyclic subgroup $\langle s_1 \rangle$ is $deg = o(s_1)(2 - 1) = o(s_1) = 2x$ and the degree of vertices in the vertex set corresponding to cosets of the cyclic subgroup $\langle s_2 \rangle$ is $deg = o(s_2)(2 - 1) = o(s_1) = 2y$. Thus, all of the vertices of $\Gamma(G, S)$ have even degree. Hence $\Gamma(G, S)$ is Eulerian. □

**Proposition 2.4.9.** *[10] The G-graph $\Gamma(G, S)$ of the group $G$ with the generating set $S$ is connected.*

As such, all the graphs under consideration in this thesis are connected.

The following result also appears in [34].

**Theorem 2.4.10.** *For any subgroup $H$ of $G$, $\Gamma(H, S_H)$ is a subgraph of $\Gamma(G, S)$ where $S_H \subseteq S$.*

*Proof.* Let $S_H \subseteq S$, and $H$ be a subgroup of $G$. Suppose $S_H = \{s_{h_1}, s_{h_2}, ..., s_{h_k}\} \subseteq S = \{s_1, s_2, ..., s_n\}$ and the group $H = \{e, h_1, ..., h_j\} \subseteq G = \{e, g_1, ..., g_n\}$.

The cosets $\langle s_{h_i} \rangle h_i, h_i \in H$ of the cyclic subgroup $\langle s_{h_i} \rangle$ have cardinality $o(s_{h_i})$. These cosets correspond to vertices in $\Gamma(H, S_H)$. Since $s_{h_i} \in S$ and $h_i \in G$, the cosets $\langle s_{h_i} \rangle h_i \in \{\langle s_i \rangle g_i\}$, which are cosets (vertices) of $\Gamma(G, S)$. The subgraph of $\Gamma(G, S)$ induced by the vertices (cosets) $\langle s_{h_i} \rangle h_i$ is $\Gamma(H, S_H)$ since, $\{\langle s_{h_i} \rangle h_i, \langle s_{h_j} \rangle h_j\}$ is a $k$-edge if and only if $|\langle s_{h_i} \rangle h_i \cap \langle s_{h_j} \rangle h_j| = k$. □

# Chapter 3

# Dependence on Generating set

In this chapter, we highlight the results which demonstrate how the $G$-graph depends on the generating set.

We had in proposition 2.4.1 that if $|S| = k$, then the G-Graph $\Gamma(G, S)$ is $k$-partite. Corollary 3.0.11 is immediate.

**Corollary 3.0.11.** *If $S'$ and $S''$ are both generating sets for the group $G$ where $|S'| = m$ and $|S''| = n$ with $m \neq n$, then the G-Graphs $\Gamma(G, S')$ and $\Gamma(G, S'')$ are not isomorphic. We write as $\Gamma(G, S') \not\cong \Gamma(G, S'')$*

*Proof.* Clearly, $\Gamma(G, S')$ will yield an $m$-partite $G$-graph while $\Gamma(G, S'')$ yields an $n$-partite $G$-graph. □

**Theorem 3.0.12.** *[10] Let $G_1$ be a group generated by $S'$ and $G_2$ be a group generated by $S''$. If $(G_1, S') \cong (G_2, S'')$ then $\Gamma(G_1, S') \cong \Gamma(G_2, S'')$ where the first $\cong$ is a group isomorphism and the second $\cong$ is a graph isomorphism.*

**Theorem 3.0.13.** *Suppose $S'$ and $S''$ are both generating sets for the group $G$, with $|S'| = |S''| = m$. If $S'$ and $S''$ have elements of the same order, then the G-graphs $\Gamma(G, S')$ and $\Gamma(G, S'')$ are isomorphic (that is, $\Gamma(G, S') \cong \Gamma(G, S'')$).*

*Proof.* Suppose $S'$ and $S''$ have elements of the same order then there is a group isomorphism $\varphi : (G, S') \to (G, S'')$ with $\varphi(s_i') = s_j'' \Leftrightarrow o(s_i') = o(s_j'')$ $s_i' \in S', s_j'' \in S''$. Hence by Theorem 3.0.12, $\Gamma(G, S') \cong \Gamma(G, S'')$. □

We consider some examples to illustrate theorem 3.0.13.

**Example 3.0.14.** If we consider the cyclic group of order 6, $G = \{e, a, a^2, a^3, a^4, a^5\}$, with generating set $S' = \{a\}$, or $\{a^5\}$. Clearly, the $G$-graph $\Gamma(G, S')$ is a vertex since it has only one coset of the the cyclic subgroup $\langle a \rangle$. $G$ can also be generated by $S'' = \{a^2, a^3\}$ with $o(a^2) = 3$ and $o(a^3) = 2$. The $G$-graph $\Gamma(G, S'')$ is isomorphic to $K_{2,3}$. Moreover, $G$ can also be generated by $S''' = \{a^3, a^4\}$ with $o(a^3) = 2$ and $o(a^4) = 3$. The $G$-graph $\Gamma(G, S''')$ is also isomorphic to $K_{2,3}$.

**Example 3.0.15.** The group $G = \mathbb{Z}_3 \times \mathbb{Z}_3$ has generating sets $S' = \{(1,0), (0,1)\}$, $S'' = \{(1,1), (1,0)\}$ and $S''' = \{(1,1), (0,1)\}$. $|S'| = |S''| = |S'''|$ and the elements in the generating sets have the same order 3. Hence $\Gamma(G, S') \cong \Gamma(G, S'') \cong \Gamma(G, S''') \cong K_{3,3}$. Fig. 3.1 is the $G$-graph $\Gamma(\mathbb{Z}_3 \times \mathbb{Z}_3, \{(1,1), (0,1)\})$.

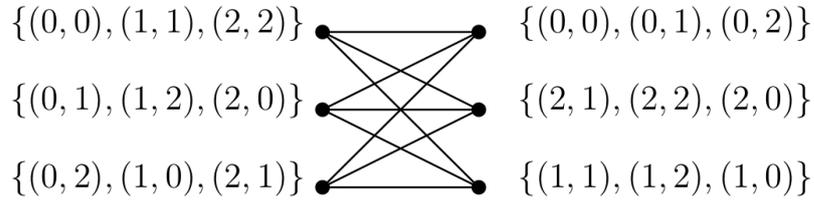

**Figure 3.1:** $\Gamma(\mathbb{Z}_3 \times \mathbb{Z}_3, \{(1,1), (0,1)\})$

**Example 3.0.16.** The Klein group $V$, is the product of two cyclic groups of order 2. $V = \{e, a, b, ab\}$ such that $ab = ba$ and $o(a) = 2 = o(b) = o(ab)$. The set $S = \{a, b, ab\}$ generates $V$. It is shown in [10] that, the $G$-graph $\Gamma(V, S)$ is isomorphic to the octahedral graph which is a 3-partite symmetric quartic graph. It is easy to see that the set $S'' = \{a, b\}$ also generates $V$. Let us construct $\Gamma(V, S'')$.

The cosets of the cyclic subgroup $\langle a \rangle$ are

$$\langle a \rangle e = \{e, a\} \text{ and } \langle a \rangle b = \{b, ab\}.$$

The cosets of the cyclic subgroup $\langle b \rangle$ are

$$\langle b \rangle e = \{e, b\} \text{ and } \langle b \rangle a = \{a, ba\} = \{a, ab\}.$$

The $G$-graph $\Gamma(V, S'')$ is isomorphic to the square, a bipartite symmetric 2-regular graph (that is, $K_{2,2}$). See Fig. 3.2.

**Example 3.0.17.** Now we consider the dihedral group $D_{2n}$. Let $r$ be a rotation through an angle of $\frac{2\pi}{n}$ and $s$ be a reflection in one particular axis. Then, $S = \{r, s\}$

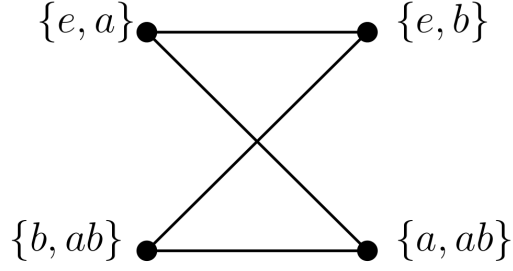

**Figure 3.2:** $\Gamma(V, S'')$

generates $D_{2n}$ where $o(r) = n$ and $o(s) = 2$. Thus, $D_{2n}$ has the presentation

$$\langle r, s | r^n = s^2 = e, srs = r^{-1} \rangle.$$

We note that $srs = r^{-1} \Rightarrow srs = r^{n-1} \Rightarrow sr = r^{n-1}s$. It is shown in [10] that the $G$-graph $\Gamma(D_{2n}, S)$ of $D_{2n}$, is isomorphic to $K_{2,n}$.

Now note that $D_{2n}$ can also be generated by two different reflections, say $s$ and $t$, that is, $S'' = \{s, t\}$. Then $D_{2n}$ has the presentation

$$\langle s, t | s^2 = t^2 = e, (ts)^n = (st)^n = e \rangle. \tag{3.1}$$

We show below that the $G$-graph $\Gamma(D_{2n}, S'')$ is isomorphic to a bipartite $C_{2n}$ graph.

We note that the presentation 3.1 implies the following relations;

1. $s(ts)^i = (ts)^{n-i}s$
2. $(st) = (ts)^{n-1}$
3. $t(ts)^i = (ts)^{n-i+1}s$

The cosets of the cyclic subgroup $\langle s \rangle$ are

$$
\begin{aligned}
\langle s \rangle e &= \{e, s\} \\
\langle s \rangle ts &= \{ts, sts\} = \{ts, (ts)^{n-1}s\} \\
\langle s \rangle (ts)^2 &= \{(ts)^2, s(ts)^2\} = \{(ts)^2, (ts)^{n-2}s\} \\
\langle s \rangle (ts)^3 &= \{(ts)^3, s(ts)^3\} = \{(ts)^3, (ts)^{n-3}s\} \\
&\vdots \\
\langle s \rangle (ts)^{n-2} &= \{(ts)^{n-2}, s(ts)^{n-2}\} = \{(ts)^{n-2}, (ts)^2 s\} \\
\langle s \rangle (ts)^{n-1} &= \{(ts)^{n-1}, s(ts)^{n-1}\} = \{(ts)^{n-1}, (ts)s\} = \{(ts)^{n-1}, t\}.
\end{aligned}
$$

So there are $n$ cosets of the cyclic subgroup $\langle s \rangle$.

The cosets of the cyclic subgroup $\langle t \rangle$ are

$$\begin{aligned}
\langle t \rangle\, e &= \{e, t\} \\
\langle t \rangle\, ts &= \{ts, t(ts)\} = \{ts, s\} \\
\langle t \rangle\, (ts)^2 &= \{(ts)^2, t(ts)^2\} = \{(ts)^2, (ts)^{n-1}s\} \\
\langle t \rangle\, (ts)^3 &= \{(ts)^3, t(ts)^3\} = \{(ts)^3, (ts)^{n-2}s\} \\
&\vdots \\
\langle t \rangle\, (ts)^{n-2} &= \{(ts)^{n-2}, t(ts)^{n-2}\} = \{(ts)^{n-2}, (ts)^3 s\} \\
\langle t \rangle\, (ts)^{n-1} &= \{(ts)^{n-1}, t(ts)^{n-1}\} = \{(ts)^{n-1}, (ts)^2 s\}.
\end{aligned}$$

So there are $n$ cosets of the cyclic subgroup $\langle t \rangle$.

Since $|S| = 2$, the graph is bipartite. All $n$ cosets of the cyclic subgroup $\langle s \rangle$ are of the form
$$\langle s \rangle\, (ts)^{n-i} = \{(ts)^{n-i}, (ts)^i s\}, \text{ where } 0 \leq i \leq n-1.$$

And all $n$ cosets of the cyclic subgroup $\langle t \rangle$ are of the form
$$\langle t \rangle\, (ts)^{n-i} = \{(ts)^{n-i}, (ts)^{n-i+1} s\}, \text{ where } 0 \leq i \leq n-1.$$

Now
$$\langle s \rangle\, (ts)^{n-i} \cap \langle t \rangle\, (ts)^{n-i} = (ts)^{n-i}$$

and
$$\langle s \rangle\, (ts)^{n-i} \cap \langle t \rangle\, (ts)^{n-i+1} = (ts)^i s.$$

Thus every vertex $(s)(ts)^{n-i}$ has degree 2.

Similarly, every vertex $(t)(ts)^{n-i}$ has degree 2 since
$$\langle t \rangle\, (ts)^{n-i} \cap \langle s \rangle\, (ts)^{n-i} = (ts)^{n-i}$$

and
$$\langle t \rangle\, (ts)^{n-i} \cap \langle s \rangle\, (ts)^{n-i-1} = (ts)^{i+1} s.$$

It can easily be verified that a graph on an even number of vertices, each of degree 2 is a cycle graph. Thus $\Gamma(D_{2n}, S)$, where $S = \{s, t\}$, is isomorphic to $C_{2n}$.

This leads us to the following proposition.

**Proposition 3.0.18.** *Let $D_{2n}$ be the dihedral group generated by two different reflections, say $s$ and $t$, that is, $S = \{s,t\}$. Then $D_{2n}$ has the presentation*

$$\langle s,t | s^2 = t^2 = e, (ts)^n = (st)^n = e \rangle.$$

*The G-graph $\Gamma(D_{2n}, S)$ is isomorphic to a bipartite $C_{2n}$ graph.*

Fig. 3.3 is an example, when $n = 4$.

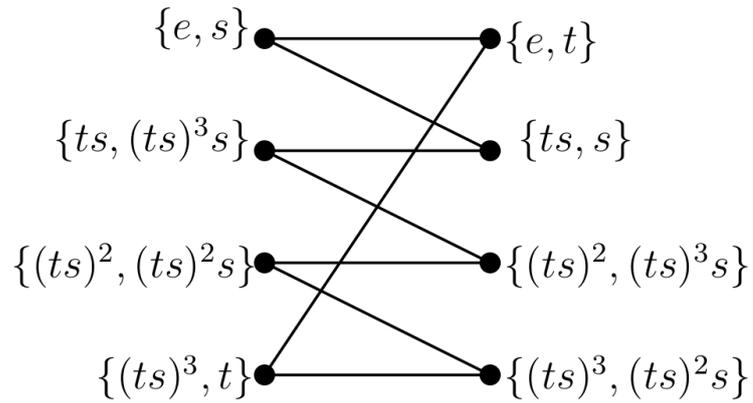

**Figure 3.3:** $\Gamma(D_8, S)$

# Chapter 4

# The Symmetric and Alternating Groups

## 4.1 Introduction

We shall consider the *symmetric group* $S_n$, the group of permutations of the set $X = \{1, 2, 3, ..., n\}$.

The *alternating group* $A_n$ is a subgroup of the symmetric group $S_n$. It consists of the elements in $S_n$ that can be written as an even number of transpositions. The order of $S_n$ is $n!$ and the order of $A_n$ is $\frac{n!}{2}$.

When a group has a single element generating set, we have the trivial $G$-graph. Hence, in all cases, we shall consider the symmetric group $S_n$, where $n \geq 3$.

The symmetric and alternating groups have several generating sets, and so we consider the $G$-graph of the symmetric and the alternating groups generated by the various generating sets.

## 4.2 The $G$-graph of $S_n$ generated by all its transpositions

In this section, we consider the symmetric group $S_n$ generated by its transpositions.

**Theorem 4.2.1.** *[17] For $n \geq 2$, $S_n$ is generated by its transpositions.*

We note that, there are $\binom{n}{2} = \frac{n(n-1)}{2}$ transpositions, i.e., $|S| = \frac{n(n-1)}{2}$. Thus, by construction, the $G$-graph is $\frac{n(n-1)}{2}$- partite. The transpositions have order 2.

Therefore, the following hold:

1. There are $\frac{n(n-1)}{2}$ cyclic subgroups $\langle s_i \rangle$.

2. The number of distinct cosets in any cyclic subgroup $\langle s_i \rangle$ is $\frac{|S_n|}{o(s_i)} = \frac{n!}{2}$.

3. All $\frac{n!}{2}$ cosets contain exactly 2 distinct elements of the symmetric group.

4. Each cyclic subgroup $\langle s_i \rangle$ has unique cosets.

**Theorem 4.2.2.** *The $G$-graph of $S_n$, with $S$, the generating set consisting of all the transpositions, is isomorphic to the $\frac{n(n-1)}{2}$- partite, $(n(n-1)-2)$-regular graph on $\frac{n!}{2} \times \frac{n(n-1)}{2}$ vertices, with $\frac{n(n-1)n!}{2}$ edges.*

*Proof.* The cardinality of the generating set is $|S| = \frac{n(n-1)}{2}$, hence it is clear the $G$-graph is $\frac{n(n-1)}{2}$-partite. Now we show that, it is an $(n(n-1)-2)$-regular graph. The transpositions have order 2, thus $o(s_i) = 2$. Let $v_i = \langle s_i \rangle g$.

$$\begin{aligned} deg(v_i) &= o(s)(|S|-1) \\ deg(v_i) &= 2(|S|-1) \\ &= 2(\tfrac{n(n-1)}{2} - 1) \\ deg(v_i) &= n(n-1) - 2. \end{aligned}$$

There are $\frac{n!}{2}$ distinct cosets in all $\frac{n(n-1)}{2}$ cyclic subgroups $\langle s_i \rangle$. Since each distinct coset forms a vertex, there shall be $\frac{n!}{2} \times \frac{n(n-1)}{2}$ vertices in all. The number of edges from Corollary 2.4.4 is $\frac{n(n-1)}{2}|G| = \frac{n(n-1)n!}{2}$. □

**Example 4.2.3.** The $G$-graph $\Gamma(S_3, S)$ of $S_3$, where $S$ is the set of transpositions, i.e., $S = \{(12), (13), (23)\}$, is isomorphic to a 3-partite, 4-regular graph on 9 vertices with 18 edges.

We shall show this constructively.

The elements of $S_3$ are $S_3 = \{(e), (12), (13), (23), (123), (132)\}$. The generating set is $S = \{(12), (13), (23)\}$.

There are 3 cosets $\langle (12) \rangle g, g \in S_3$. They are

$$\langle (12) \rangle e = \{e, (12)\}$$
$$\langle (12) \rangle (13) = \{(13), (132)\}$$
$$\langle (12) \rangle (23) = \{(23), (123)\}.$$

There are 3 cosets $\langle (13) \rangle g, g \in S_3$. They are

$$\langle (13) \rangle e = \{e, (13)\}$$
$$\langle (13) \rangle (12) = \{(12), (123)\}$$
$$\langle (13) \rangle (23) = \{(23), (132)\}.$$

And finally, there are 3 cosets $\langle (23) \rangle g, g \in S_3$. They are

$$\langle (23) \rangle e = \{e, (23)\}$$
$$\langle (23) \rangle (12) = \{(12), (132)\}$$
$$(\langle (23) \rangle) (13) = \{(13), (123)\}.$$

The graph is as shown in Fig. 4.1.

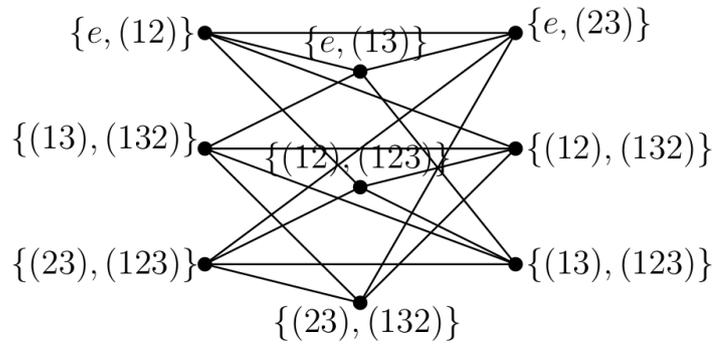

**Figure 4.1:** $\Gamma(S_3, S)$

**Example 4.2.4.** The $G$-graph $\Gamma(S_4, S)$ of $S_4$, where $S$ is the set of transpositions, i.e., $S = \{(12), (13), (14), (23), (24), (34)\}$, is isomorphic to a 6-partite, 10-regular graph on 72 vertices.

Each of the 6 vertex sets contain 12 vertices, each vertex is adjacent to two vertices in the other 5 vertex sets.

**Example 4.2.5.** The $G$-graph $\Gamma(S_5, S)$ of $S_5$, where $S$ is the set of transpositions, i.e., $S = \{(12),(13),(14),(15),(23),(24),(25),(34),(35),(45)\}$, is isomorphic to a 10-partite, 18-regular graph on 600 vertices.

Each of the 10 vertex sets contain 60 vertices, each vertex is adjacent to two vertices in the other 9 vertex sets.

## 4.3 The $G$-graph of $S_n$ generated by $(n-1)$ transpositions

In this section, we consider the symmetric group $S_n$ generated by some of its transpositions.

**Theorem 4.3.1.** *[17] For $n \geq 2$, $S_n$ is generated by the $(n-1)$ transpositions*

$$(12),(13),(14),...,(1n).$$

**Theorem 4.3.2.** *[17] For $n \geq 2$, $S_n$ is generated by the $(n-1)$ transpositions*

$$(12),(23),(34),...,(n-1\,n).$$

**Theorem 4.3.3.** *The $G$-graph $\Gamma(S_n, S)$ of $S_n$, with $S$, the generating set consisting of $(n-1)$ transpositions, i.e., $S = \{(12),(23),(34),...,(n-1\,n)\}$ is isomorphic to the $(n-1)$- partite, $(2n-4)$-regular graph on $\frac{n(n-1)(n-1)!}{2}$ vertices.*

*Proof.* Since the cardinality of $S$ is $n-1$, the graph is $(n-1)$- partite. Moreover, the elements of the generating set have the same order 2, thus the graph $\Gamma(S_n, S)$ is regular and for all $v \in V$ $deg(v) = 2((n-1)-1) = 2(n-2) = 2n-4$. There are $\frac{n!}{2}$ distinct cosets in all $(n-1)$ cyclic subgroups $\langle s_i \rangle$. Since each distinct coset forms a vertex, there shall be $\frac{n!}{2} \times (n-1) = \frac{n(n-1)(n-1)!}{2}$ vertices in all. $\square$

**Corollary 4.3.4.** *If $S = \{(12),(23),(34),...,(n-1n)\}$ and $S^* = (12),(13),(14),...,(1n)$, then*

$$\Gamma(S_n, S) \cong \Gamma(S_n, S^*).$$

*Proof.* $S^*$ generates $S_n$ by theorem 4.3.1 and $S$ generates $S_n$ by theorem 4.3.2. $|S^*| = |S| = n-1$. Each $s \in S$ has order 2 and each $s^* \in S^*$ has order 2, thus $(S_n, S) \cong (S_n, S^*)$ and so by theorem 3.0.13, $\Gamma(S_n, S) \cong \Gamma(S_n, S^*)$. □

**Example 4.3.5.** The $G$-graph $\Gamma(S_3, S)$ of $S_3$, where $S = \{(12), (23)\}$, is isomorphic to a bipartite, 2-regular graph on 6 vertices, i.e. $C_6$.

We shall explicitly construct $\Gamma(S_3, S)$.

The elements of $S_3$ are $S_3 = \{(e), (12), (13), (23), (123), (132)\}$.

The generating set is $S = \{(12), (23)\}$.

There are 3 cosets $\langle (12) \rangle g, g \in S_3$. They are

$$\langle (12) \rangle e = \{e, (12)\}$$
$$\langle (12) \rangle (13) = \{(13), (132)\}$$
$$\langle (12) \rangle (23) = \{(23), (123)\}.$$

And finally, there are 3 cosets $\langle (23) \rangle g, g \in S_3$. They are

$$\langle (23) \rangle e = \{e, (23)\}$$
$$\langle (23) \rangle (12) = \{(12), (132)\}$$
$$\langle (23) \rangle (13) = \{(13), (123)\}.$$

The graph is as shown in Fig. 4.2.

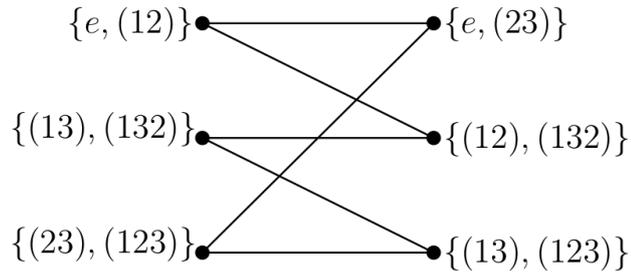

**Figure 4.2:** $\Gamma(S_3, S) \cong C_6$

**Example 4.3.6.** The $G$-graph $\Gamma(S_4, S)$ of $S_4$, where $S = \{(12), (23), (34)\}$, is isomorphic to a 3-partite, 4-regular graph on 36 vertices.

Each of the 3 vertex sets contain 12 vertices and each vertex is adjacent to two vertices in the other 2 vertex sets.

**Example 4.3.7.** The $G$-graph $\Gamma(S_5, S)$ of $S_5$, where $S = \{(12), (23), (34), (45)\}$, is isomorphic to a 4-partite, 6-regular graph on 240 vertices.

Each of the 4 vertex sets contain 60 vertices and each vertex is adjacent to two vertices in each of the other 3 vertex sets.

## 4.4 The $G$-graph of $S_n$ generated by the transposition $(12)$ and the $n$-cycle $(12...n)$

In this section, we consider the symmetric group $S_n$ generated by two elements.

**Theorem 4.4.1.** *[17] For $n \geq 2$, $S_n$ is generated by the transposition $(12)$ and the $n$-cycle $(12...n)$.*

**Theorem 4.4.2.** *The $G$-graph $\Gamma(S_n, S)$ of $S_n \forall n > 2$, with $S = \{(12), (12...n)\}$ is isomorphic to a semi-regular bipartite graph with vertex set $V = (V_1 \cup V_2)$, $|V_1| = \frac{n!}{2}$, $|V_2| = (n-1)!$ and $deg(v_{1_i}) = 2 \ \forall v_{1_i} \in V_1$, $deg(v_{2_j}) = n \ \forall v_{2_j} \in V_2$.*

*Proof.* Since $|S| = 2$, the graph is bipartite. Moreover, the two elements $(12), (12...n)$ of the generating set $S$ have orders 2 and $n$ respectively. Thus there are $\frac{n!}{2}$ distinct cosets in the cyclic subgroup $\langle (12) \rangle$ and since each distinct coset forms a vertex, we have $|V_1| = \frac{n!}{2}$. Similarly, there are $\frac{n!}{n}$ distinct cosets in the cyclic subgroup $\langle (12...n) \rangle$ and since each distinct coset forms a vertex, we have $|V_2| = \frac{n!}{n} = (n-1)!$.
Now, $\forall v_{1_i} \in V_1$, $deg(v_{1_i}) = o((12))(|S| - 1) = 2(2 - 1) = 2$ and $\forall v_{2_j} \in V_2$, $deg(v_{2_j}) = o((12...n))(|S| - 1) = n(2 - 1) = n$. $\square$

**Remark 4.4.3.** We note that, the number of edges in this graph is $n!$, the order of $S_n$. So by construction, each edge corresponds to an element in $S_n$.

**Example 4.4.4.** The $G$-graph $\Gamma(S_3, S)$ of $S_3$, where $S = \{(12), (123)\}$, is isomorphic to the complete bipartite graph $K_{23}$.

We shall explicitly construct $\Gamma(S_3, S)$.

The group $S_3 = \{(e), (12), (13), (23), (123), (132)\}$.

The generating set is $S = \{(12), (123)\}$, the order of $(12)$ is $2$ and the order of $(123)$ is $3$. The cardinality of $S$ is $2$.

There are $3$ cosets $\langle (12) \rangle\, g, g \in S_3$. They are

$$\langle (12) \rangle\, e = \{e, (12)\}$$
$$\langle (12) \rangle\, (13) = \{(13), (132)\}$$
$$\langle (12) \rangle\, (23) = \{(23), (123)\}.$$

And finally, there are $2$ cosets $\langle (123) \rangle\, g, g \in S_3$. They are

$$\langle (123) \rangle\, e = \{e, (123), (132)\}$$
$$\langle (123) \rangle\, (12) = \{(12), (13), (23)\}.$$

The graph is as shown in Fig. 4.3.

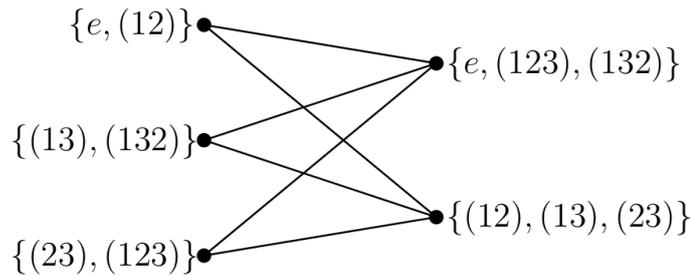

**Figure 4.3:** $\Gamma(S_3, S) \cong K_{23}$

**Example 4.4.5.** The $G$-graph $\Gamma(S_4, S)$ of $S_4$, where $S = \{(12), (1234)\}$, is isomorphic to a semi-regular bipartite graph with $V = (V_1 \cup V_2)$ where $|V_1| = 12$, $|V_2| = 6$, $\forall a \in V_1, deg(a) = 2$, $\forall b \in V_2, deg(b) = 4$ and $|E| = 24$.

**Example 4.4.6.** The $G$-graph $\Gamma(S_5, S)$ of $S_5$, where $S = \{(12), (12345)\}$, is isomorphic to a semi-regular bipartite graph with $V = (V_1 \cup V_2)$ where $|V_1| = 60$, $|V_2| = 24$, $\forall a \in V_1, deg(a) = 2$, $\forall b \in V_2, deg(b) = 5$ and $|E| = 120$.

## 4.5 The $G$-graph of $S_n$ generated by the transposition $(12)$ and the $(n-1)$-cycle $(23...n)$

In this section, we also consider the symmetric group $S_n$ generated by two elements.

**Theorem 4.5.1.** *[17] For $n \geq 2$, $S_n$ is generated by the transposition $(12)$ and the $(n-1)$-cycle $(23...n)$.*

**Theorem 4.5.2.** *For all $n > 2$, the $G$-graph $\Gamma(S_n, S)$ of $S_n$, with $S = \{(12), (23...n)\}$ is isomorphic to a semi-regular bipartite graph with vertex set $V = (V_1 \cup V_2)$, $|V_1| = \frac{n!}{2}$, $|V_2| = \frac{n!}{(n-1)}$ and $deg(v_{1_i}) = 2\ \forall v_{1_i} \in V_1$, $deg(v_{2_j}) = n - 1\ \forall v_{2_j} \in V_2$.*

*Proof.* Since $|S| = 2$, the graph is bipartite. Moreover, the two elements $(12), (23...n)$ of the generating set $S$ have orders $2$ and $n - 1$ respectively. Thus there are $\frac{n!}{2}$ distinct cosets in the cyclic subgroup $\langle (12) \rangle$ and since each distinct coset forms a vertex, we have $|V_1| = \frac{n!}{2}$. Similarly, there are $\frac{n!}{n-1}$ distinct cosets in the cyclic subgroup $\langle (23...n) \rangle$ and since each distinct coset forms a vertex, we have $|V_2| = \frac{n!}{n-1}$. Now, $\forall v_{1_i} \in V_1$, $deg(v_{1_i}) = o((12))(|S| - 1) = 2(2 - 1) = 2$ and $\forall v_{2_j} \in V_2$ $deg(v_{2_j}) = o((23...n))(|S| - 1) = (n - 1)(2 - 1) = n - 1$. □

**Remark 4.5.3.** This graph is regular only when $n = 3$ in which case $deg(v_{2_j}) = n - 1$ becomes $deg(v_{2_j}) = 3 - 1 = 2$.

**Remark 4.5.4.** We note that, the number of edges in this graph is $n!$, the order of $S_n$. So by construction, each edge corresponds to an element in $S_n$.

**Example 4.5.5.** The $G$-graph $\Gamma(S_3, S)$ of $S_3$, where $S = \{(12), (23)\}$, is isomorphic to the complete bipartite graph $C_6$.

This is Example 4.3.5.

**Example 4.5.6.** The $G$-graph $\Gamma(S_4, S)$ of $S_4$, where $S = \{(12), (234)\}$, is isomorphic to a semi-regular bipartite graph with $V = (V_1 \cup V_2)$ where $|V_1| = 12$, $|V_2| = 8$, $\forall a \in V_1, deg(a) = 2$, $\forall b \in V_2, deg(b) = 3$ and $|E| = 24$.

We shall explicitly construct $\Gamma(S_4, S)$.

There are 12 cosets $\langle (12) \rangle g, g \in S_4$. They are

$$\langle (12) \rangle e = \{e, (12)\}$$
$$\langle (12) \rangle (13) = \{(13), (132)\}$$
$$\langle (12) \rangle (14) = \{(14), (142)\}$$
$$\langle (12) \rangle (23) = \{(23), (123)\}$$
$$\langle (12) \rangle (24) = \{(24), (124)\}$$
$$\langle (12) \rangle (34) = \{(34), (12)(34)\}$$
$$\langle (12) \rangle (134) = \{(134), (1342)\}$$
$$\langle (12) \rangle (143) = \{(143), (1432)\}$$
$$\langle (12) \rangle (234) = \{(234), (1234)\}$$
$$\langle (12) \rangle (243) = \{(243), (1243)\}$$
$$\langle (12) \rangle (1324) = \{(1324), (13)(24)\}$$
$$\langle (12) \rangle (1423) = \{(1423), (14)(23)\}.$$

There are 8 cosets $\langle (234) \rangle g, g \in S_4$. They are

$$\langle (234) \rangle e = \{e, (234), (243)\}$$
$$\langle (234) \rangle (13) = \{(13), (1423)), (1243)\}$$
$$\langle (234) \rangle (14) = \{(14), (1234), (1324)\}$$
$$\langle (234) \rangle (23) = \{(23), (24), (34)\}$$
$$\langle (234) \rangle (12) = \{(12), (1342), (1432)\}$$
$$\langle (234) \rangle (123) = \{(123), (13)(24), (143)\}$$
$$\langle (234) \rangle (132) = \{(132), (142), (12)(34)\}$$
$$\langle (234) \rangle (124) = \{(124), (134), (14)(23)\}.$$

The graph is as shown in Fig. 4.4.

**Example 4.5.7.** The $G$-graph $\Gamma(S_5, S)$ of $S_5$, where $S = \{(12), (2345)\}$, is isomorphic to a semi-regular bipartite graph with $V = (V_1 \cup V_2)$ where $|V_1| = 60$, $|V_2| = 30$, $\forall a \in V_1, deg(a) = 2$, $\forall b \in V_2, deg(b) = 4$ and $|E| = 120$.

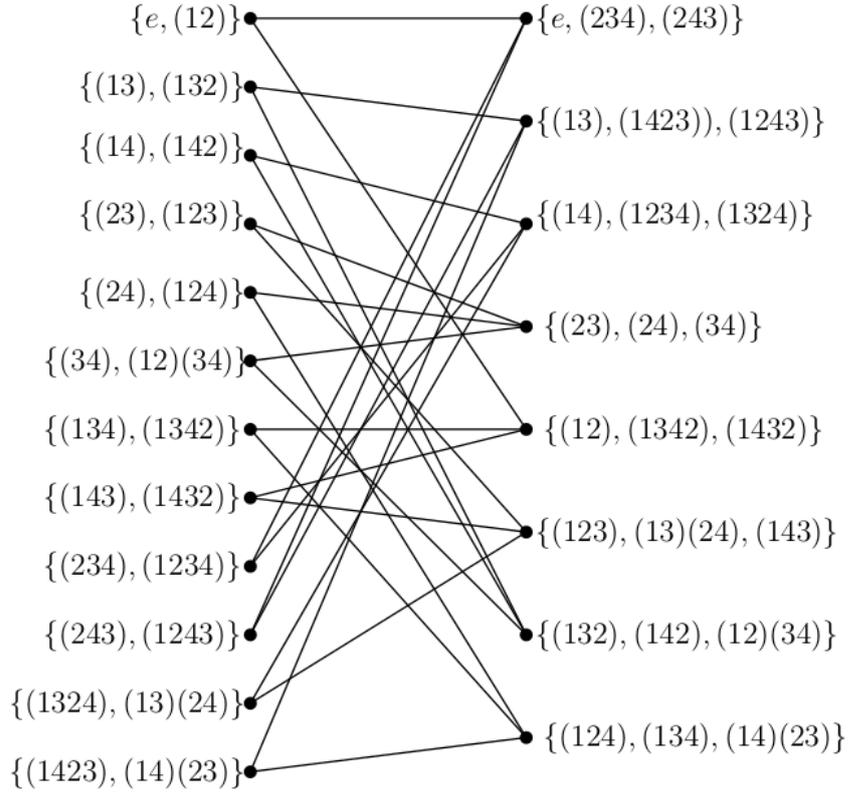

**Figure 4.4:** $\Gamma(S_4, S)$

Now we shift our attention to the alternating group $A_n$. It is known [17] that $A_n$ can be generated by all of its $\frac{n(n-1)(n-2)}{3}$ 3-cycles but we want to consider generating sets with fewer elements.

## 4.6 The $G$-graph of the Alternating group generated by $(12i)$'s

In this section, we consider the Alternating group $A_n$ generated by some of its 3-cycles.

**Theorem 4.6.1.** *[17] For $n \geq 3$, the group $A_n$ is generated by 3-cycles of the form $(12i)$.*

**Theorem 4.6.2.** *For all $n > 3$, the $G$-graph of $A_n$ generated by 3-cycles of the form $(12i)$ is isomorphic to a $(n-2)$-partite, $3(n-3)$-regular graph such that, $|E| = \frac{n!(n-2)(n-3)}{4}$ and $|V| = \frac{n!(n-2)}{6}$.*

*Proof.* There are $n-2$, 3-cycles of the form $(12i)$ in $A_n$. $o(12i) = 3$ and so the degree of each vertex is $3((n-2)-1) = 3(n-3)$. By corollary 2.4.4, $|E| = \frac{(n-2)(n-3)}{2} \times |A_n| = \frac{n!(n-2)(n-3)}{4}$ and by corollary 2.4.3, $|V| = \sum_{i=1}^{n-2} \frac{|G|}{o(s_i)} = \frac{n!}{6}(n-2) = \frac{n!(n-2)}{6}$ since $o(s_i) = 3$ for all $s_i \in S$. $\square$

**Example 4.6.3.** The $G$-graph $\Gamma(A_4, S)$ of $A_4$, where $S = \{(123), (124)\}$, is isomorphic to a bipartite, 3-regular graph on 8 vertices with 12 edges.

We shall explicitly construct $\Gamma(A_4, S)$.

The elements of $A_4$ are $A_4 = \{(e), (12)(34), (13)(24), (14)(23), (123), (132), (124),$
$(142), (134), (143), (234), (243)\}$. The generating set is $S = \{(123), (124)\}$.

There are 4 cosets $\langle (123) \rangle g, g \in A_4$. They are

$$\langle (123) \rangle e = \{e, (123), (132)\}$$
$$\langle (123) \rangle (124) = \{(124), (13)(24), (243)\}$$
$$\langle (123) \rangle (142) = \{(142), (143), (14)(23)\}$$
$$\langle (123) \rangle (134) = \{(134), (234), (12)(34)\}.$$

And finally, there are 4 cosets $\langle (124) \rangle g, g \in A_4$. They are

$$\langle (124) \rangle e = \{e, (124), (142)\}$$
$$\langle (124) \rangle (123) = \{(123), (14)(23), (234)\}$$
$$\langle (124) \rangle (132) = \{(132), (134), (13)(24)\}$$
$$\langle (124) \rangle (243) = \{(243), (12)(34), (143)\}.$$

The graph is as shown in Fig. 4.5.

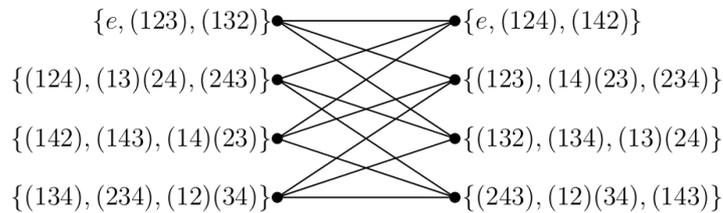

**Figure 4.5:** $\Gamma(A_4, S)$

## 4.7 The $G$-graph of the Alternating group generated by $(i \quad i+1 \quad i+2)$'s

In this section, we consider the Alternating group $A_n$ generated by another set of its 3-cycles.

**Theorem 4.7.1.** *[17] For $n \geq 3$, the consecutive 3-cycles $(i \quad i+1 \quad i+2)$, with $1 \leq i \leq n-2$, generate $A_n$.*

**Theorem 4.7.2.** *For all $n \geq 3$, the $G$-graph of $A_n$ generated by 3-cycles of the form $(i \quad i+1 \quad i+2)$ is isomorphic to a $(n-2)$-partite, $3(n-3)$-regular graph such that, $|E| = \frac{n!(n-2)(n-3)}{4}$ and $|V| = \frac{n!(n-2)}{6}$.*

*Proof.* There are $n-2$, 3-cycles of the form $(i \quad i+1 \quad i+2)$ in $A_n$. $o(i \quad i+1 \quad i+2) = 3$ and so the degree of each vertex is $3((n-2)-1) = 3(n-3)$. By corollary 2.4.4 $|E| = \frac{(n-2)(n-3)}{2} \times |A_n| = \frac{n!(n-2)(n-3)}{4}$ and by corollary 2.4.3 $|V| = \sum_{i=1}^{n-2} \frac{|G|}{o(s_i)} = \frac{n!}{6}(n-2) = \frac{n!(n-2)}{6}$ since $o(s_i) = 3$ for all $s_i \in S$. □

**Remark 4.7.3.** Theorems 4.7.2 and 4.6.2 also provide us examples of theorem 3.0.13.

**Example 4.7.4.** The $G$-graph $\Gamma(A_4, S)$ of $A_4$, where $S = \{(123), (234)\}$, is isomorphic to a bipartite, 3-regular graph on 8 vertices with 12 edges.

We shall explicitly construct $\Gamma(A_4, S)$.

The elements of $A_4$ are

$$A_4 = \{(e), (12)(34), (13)(24), (14)(23), (123), (132), (124),$$

$$(142), (134), (143), (234), (243)\}.$$

The generating set is $S = \{(123), (124)\}$.

There are 4 cosets $\langle (123) \rangle g, g \in A_4$. They are

$$\langle (123) \rangle e = \{e, (123), (132)\}$$
$$\langle (123) \rangle (124) = \{(124), (13)(24), (243)\}$$
$$\langle (123) \rangle (142) = \{(142), (143), (14)(23)\}$$
$$\langle (123) \rangle (134) = \{(134), (234), (12)(34)\}.$$

And finally, there are 4 cosets $\langle (234) \rangle g, g \in A_4$. They are

$$\langle (234) \rangle e = \{e, (234), (243)\}$$
$$\langle (234) \rangle (124) = \{(124), (134), (14)(23)\}$$
$$\langle (234) \rangle (123) = \{(123)), (13)(24), (143)\}$$
$$\langle (234) \rangle (142) = \{(142), (12)(34), (132)\}.$$

The graph is as shown in Fig. 4.6.

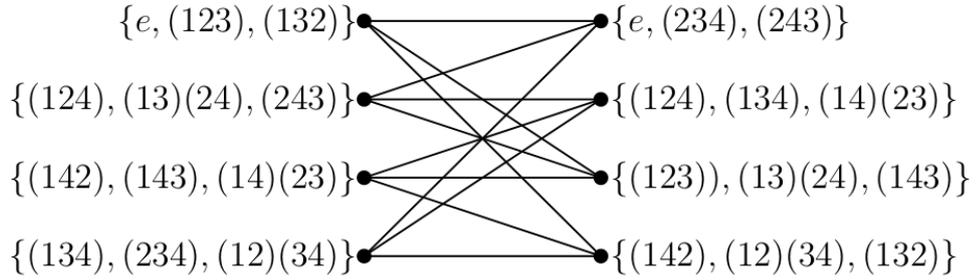

**Figure 4.6:** $\Gamma(A_4, S)$

## 4.8 The $G$-graph of the Alternating group generated by Two elements

In this section, we consider the Alternating group $A_n$ generated by two elements. This depends on whether $n$ is odd or even.

**Theorem 4.8.1.** *[17] For $n \geq 4$, $A_n$ is generated by the two elements* (123) *and*

$$\begin{cases} (12...n) & \text{if } n \text{ is odd} \\ (23...n) & \text{if } n \text{ is even.} \end{cases}$$

## 4.8.2 The $G$-graph of the Alternating group generated by $(123)$ and $(12...n)$

From theorem 4.8.1, for all $n$ odd, $A_n$ is generated by $(123)$ and $(12...n)$. The order of $(123)$ is 3 while the order of $(12...n)$ is $n$.

**Theorem 4.8.3.** *Let $n \geq 4$. For all $n$ odd, the $G$-graph $\Gamma(A_n, S)$ of $A_n$ generated by $S = \{(123), (12...n)\}$ is isomorphic to a bipartite graph with vertex set $V = V_1 \cup V_2$ where $|V_1| = \frac{n!}{6}$ and $|V_2| = \frac{n!}{2n}$ with $deg(a) = 3$ for all $a \in V_1$ and $deg(b) = n$ for all $b \in V_2$ such that, $|E| = \frac{n!}{2}$.*

*Proof.* $|S| = 2$ so $\Gamma(A_n, S)$ is bipartite and there are $\frac{n!}{2 \times 3}$ cosets $\langle(123)\rangle g \quad g \in A_n$ so $|V_1| = \frac{n!}{6}$ while there are $\frac{n!}{2 \times n}$ cosets $\langle(12...n)\rangle g \quad g \in A_n$ so $|V_2| = \frac{n!}{2n}$. For any vertex $a$ in $V_1$, $deg(a) = 3 = o(123)$ and for any vertex $b$ in $V_2$, $deg(b) = n = o(12...n)$. By corollary 2.4.4, $|E| = \frac{2(2-1)}{2} \times |A_n| = \frac{n!}{2}$. □

**Example 4.8.4.** The $G$-graph $\Gamma(A_5, S)$ of $A_5$, where $S = \{(123), (12345)\}$, is isomorphic to a bipartite graph with vertex set $V = V_1 \cup V_2$ where $|V_1| = \frac{5!}{6} = 20$ and $|V_2| = \frac{5!}{10} = 12$ with $deg(a) = 3$ for all $a \in V_1$ and $deg(b) = 5$ for all $b \in V_2$ such that, $|E| = 60$.

## 4.8.5 The $G$-graph of the Alternating group generated by $(123)$ and $(23...n)$

From theorem 4.8.1, for all $n$ even, $A_n$ is generated by $(123)$ and $(23...n)$. The order of $(123)$ is 3 while the order of $(23...n)$ is $n-1$.

**Theorem 4.8.6.** *Let $n \geq 4$. For all $n$ even, the $G$-graph of $A_n$ generated by $S = \{(123), (23...n)\}$ is isomorphic to a bipartite graph with vertex set $V = V_1 \cup V_2$ where $|V_1| = \frac{n!}{6}$ and $|V_2| = \frac{n!}{2n-2}$ with $deg(a) = 3$ for all $a \in V_1$ and $deg(b) = n - 1$ for all $b \in V_2$ such that $|E| = \frac{n!}{2}$.*

*Proof.* $|S| = 2$ so $\Gamma(A_n, S)$ is bipartite and there are $\frac{n!}{2 \times 3}$ cosets $\langle(123)\rangle g \quad g \in A_n$ so $|V_1| = \frac{n!}{6}$ while there are $\frac{n!}{2 \times (n-1)}$ cosets $\langle(23...n)\rangle g \quad g \in A_n$ so $|V_2| = \frac{n!}{2n}$. For any vertex $a$ in $V_1$, $deg(a) = 3 = o(123)$ and for any vertex $b$ in $V_2$, $deg(b) = n - 1 = o(23...n)$. Finally by corollary 2.4.4, $|E| = \frac{2(2-1)}{2} \times |A_n| = \frac{n!}{2}$. □

**Example 4.8.7.** When $n = 4$, we have precisely example 4.7.4. That is, the $G$-graph $\Gamma(A_4, S)$ of $A_4$, where $S = \{(123), (234)\}$, is isomorphic to a bipartite graph with vertex set $V = V_1 \cup V_2$ where $|V_1| = 4$ and $|V_2| = 4$ with $deg(a) = 3$ for all $a \in V_1$ and $deg(b) = 3$ for all $b \in V_2$ such that, $|E| = 12$.

**Example 4.8.8.** When $n = 6$, that is, the $G$-graph $\Gamma(A_6, S)$, where $S = \{(123), (23456)\}$, is isomorphic to a bipartite graph with vertex set $V = V_1 \cup V_2$ where $|V_1| = 120$ and $|V_2| = 72$ with $deg(a) = 3$ for all $a \in V_1$ and $deg(b) = 5$ for all $b \in V_2$ such that, $|E| = 360$.

# Chapter 5

# The Quaternions, dihedral and semi-dihedral groups

## 5.1 Introduction

In this chapter, we present the $G$-graphs of the quaternion, dihedral and semi-dihedral groups.

## 5.2 The $G$-graph of the Quaternions

The quaternion group has the presentation $Q = \langle a, b | a^4 = e, b^2 = a^2, bab^{-1} = a^{-1} \rangle$. It has order 8.

**Theorem 5.2.1.** *The $G$-graph of the quaternion group $Q$ with the presentation above, is isomorphic to the double edged complete bipartite graph $K_{22}$.*
*That is $\Gamma(Q, \{a, b\}) \cong K_{22}^2$.*

*Proof.* Since the generating set $S = \{a, b\}$ has $|S| = 2$, the graph is bipartite, say $V = V_1 \cup V_2$. Since $o(a) = 4$, $|V_1| = \frac{8}{4} = 2$ vertices, each of degree $o(a)(2-1) = 4$. Also, $o(b) = 4$, $|V_2| = \frac{8}{4} = 2$ vertices, each of degree $o(b)(2-1) = 4$. □

Alternatively, we prove theorem 5.2.1 by explicitly constructing $\Gamma(Q, \{a, b\})$.

The cosets of the cyclic subgroup $\langle a \rangle$ are the following. There are two of them.

$$\langle a \rangle e = \{e, a, a^2, a^3\}$$
$$\langle a \rangle b = \{b, ab, a^2b, a^3b\}.$$

Finally, the cosets of the cyclic subgroup $\langle a \rangle$ are the following. There are two of them.

$$\langle b \rangle e = \{e, b, a^2, a^2b\}$$
$$\langle b \rangle a = \{a, a^3b, a^3, ab\}.$$

The graph is as shown in Fig. 5.1.

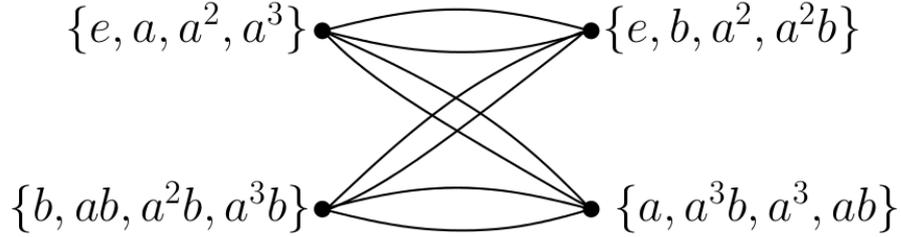

**Figure 5.1:** $\Gamma(Q, \{a,b\}) \cong K_{22}^2$

### 5.2.2 The Generalized Quaternion group

The following theorem and proof appeared in [11].

**Theorem 5.2.3.** *Let the generalized quaternion group $Q_n$ be the group with presentation*

$$\langle a, b | a^{2n} = e, b^2 = a^n, ab = ba^{2n-1} \rangle.$$

*For $S = \{a, b\}$, the graph $\Gamma(Q_n, S)$ of the generalized quaternion group is the complete double-edged bipartite graph $K_{2,n}^2$.*

## 5.3 The $G$-graph of the Dihedral group

Now we consider the dihedral group $D_{2n}$. Let $r$ be a rotation through an angle of $\frac{2\pi}{n}$ and $s$ be a reflection in one particular axis, then, $S = \{r, s\}$ generates $D_{2n}$ where

$o(r) = n$ and $o(s) = 2$. Thus $D_{2n}$ has the presentation

$$\langle r, s | r^n = s^2 = e, srs = r^{-1} \rangle.$$

We note that, $srs = r^{-1} \Rightarrow srs = r^{n-1} \Rightarrow sr = r^{n-1}s$.

The following theorem and proof appeared in [10].

**Theorem 5.3.1.** *Let $D_{2n}$ have the presentation as above. Then the $G$-graph $\Gamma(D_{2n}; S)$ of $D_{2n}$ is isomorphic to $K_{2,n}$.*

Now note that, $D_{2n}$ can also be generated by two different reflections say $s$ and $t$, i.e. $S = \{s, t\}$. Then $D_{2n}$ has the presentation

$$\langle s, t | s^2 = t^2 = e, (ts)^n = (st)^n = e \rangle.$$

We showed in Prop. 3.0.18 that, the $G$-graph $\Gamma(D_{2n}, S)$ is isomorphic to a bipartite $C_{2n}$ graph.

## 5.4 The $G$-graph of the Semi-dihedral group

The semi-dihedral group also known as the quasi-dihedral group has the presentation

$$SD_{8k} = \langle a, b | a^{4k} = b^2 = 1, \quad ba = a^{2k-1}b \rangle.$$

**Theorem 5.4.1.** *If $SD_{8k}$ is the semi-dihedral group with the presentation above, then for $S = \{a, b\}$, the $G$-graph $\Gamma(SD_{8k}; S)$ is isomorphic to a bipartite graph with $V = V_1 \cup V_2$, such that $|V_1| = 2$ and $|V_2| = 4k$ where $deg(x) = 4k \quad \forall x \in V_1$ and $deg(y) = 2 \quad \forall y \in V_2$ with $|E| = 8k$. This bipartite graph is the complete graph $K_{2,4k}$.*

*Proof.* From the presentation, the generating set is $S = \{a, b\}$ with $o(a) = 4k$ and $o(b) = 2$. Since $|S| = 2$, $\Gamma(SD_{8k}; S)$ is bipartite. Let the vertices be $V = V_1 \cup V_2$, then the distinct cosets $\langle a \rangle g \quad \forall g \in SD_{8k}$ correspond to $x_i \in V_1$ and $|V_1| = \frac{|SD_{8k}|}{o(a)} = \frac{8k}{4k} = 2$. Similarly, the distinct cosets $\langle b \rangle g \quad \forall g \in SD_{8k}$ correspond to $y_j \in V_2$ and $|V_2| = \frac{|SD_{8k}|}{o(b)} = \frac{8k}{2} = 4k$. Now for all $x_i \in V_1 \quad deg(x_i) = o(a)(|S|-1) = 4k(2-1) = 4k$ and for all $y_j \in V_1$ $deg(y) = o(b)(|S| - 1) = 2(2 - 1) = 2$. The number of edges $E$ will thus be $|E| = 2 \times 4k = 8k$. The graph with such a description is the $K_{2,4k}$ complete bipartite graph. $\square$

Now we prove theorem 5.4.1 by explicitly constructing $\Gamma(SD_{8k}; S)$.

*Proof.* The cosets of the cyclic subgroup $\langle a \rangle$ are the following. There are two of them.

$$\langle a \rangle e = \{e, a, a^2, a^3, ..., a^{4k-1}\}$$
$$\langle a \rangle b = \{b, ab, a^2b, a^3b, ..., a^{4k-1}b\}.$$

The cosets of the cyclic subgroup $\langle b \rangle$ are the following. There are $4k$ of them.

$$\langle b \rangle e = \{e, b\}$$
$$\langle b \rangle a = \{a, ba\} = \{a, a^{2k-1}b\}$$
$$\langle b \rangle a^2 = \{a^2, ba^2\} = \{a^2, a^{4k-2}b\}$$
$$\langle b \rangle a^3 = \{a^3, ba^3\} = \{a^3, a^{6k-3}b\} = \{a^3, a^{2k-3}b\}$$
$$\langle b \rangle a^4 = \{a^4, ba^4\} = \{a^4, a^{8k-4}b\} = \{a^4, a^{4k-4}b\}$$
$$\langle b \rangle a^5 = \{a^5, ba^5\} = \{a^5, a^{10k-5}b\} = \{a^5, a^{2k-5}b\}$$
$$\langle b \rangle a^6 = \{a^6, ba^6\} = \{a^6, a^{12k-6}b\} = \{a^6, a^{4k-6}b\}$$
$$\vdots$$
$$\langle b \rangle a^{4k-6} = \{a^{4k-6}, ba^{4k-6}\} = \{a^{4k-6}, a^6b\}$$
$$\langle b \rangle a^{4k-5} = \{a^{4k-5}, ba^{4k-5}\} = \{a^{4k-5}, a^{2k+5}b\}$$
$$\langle b \rangle a^{4k-4} = \{a^{4k-4}, ba^{4k-4}\} = \{a^{4k-4}, a^4b\}$$
$$\langle b \rangle a^{4k-3} = \{a^{4k-3}, ba^{4k-3}\} = \{a^{4k-3}, a^{2k+3}b\}$$
$$\langle b \rangle a^{4k-2} = \{a^{4k-2}, ba^{4k-2}\} = \{a^{4k-2}, a^2b\}$$
$$\langle b \rangle a^{4k-1} = \{a^{4k-1}, ba^{4k-1}\} = \{a^{4k-1}, a^{2k+1}b\}.$$

$|S| = 2$ so the graph is bipartite.

We must show that it is a complete graph, i.e. for all $(x, y)$, $x$ in $V_a$ and $y$ in $V_b$, there is exactly one edge between $x$ and $y$.

For all $y$ in $V_b$ there is $i$ in $\{0, 1, ..., 4k - 1\}$ such that $y = \langle b \rangle a^i$. But we have $\langle b \rangle a^i = \{a^i, a^{4k-i}b\}$ for even $i$ and $\langle b \rangle a^i = \{a^i, a^{2k-i}b\}$ for odd $i$.

If $x = \langle a \rangle b$, because $\langle a \rangle b = \{b, ab, a^2b, a^3b, ..., a^{4k-1}b\}$, we have $\langle a \rangle b \cap \langle b \rangle a^i = \{a^{4k-i}b\}$

for even $i$ and $\langle a \rangle b \cap \langle b \rangle a^i = \{a^{2k-i}b\}$ for odd $i$.

If $x = \langle a \rangle e$, because $\langle a \rangle e = \{e, a, a^2, a^3, ..., a^{4k-1}\}$, we have $\langle a \rangle e \cap \langle b \rangle a^i = \{a^i\}$.

So, for all $x$ in $V_a$ and $y$ in $V_b$, there is exactly one edge between $x$ and $y$. □

**Corollary 5.4.2.** *The $G$-graph $\Gamma(SD_{8k}; S)$ of the semi-dihedral group, is Eulerian.*

**Example 5.4.3.** For $k = 2$, the $G$-graph $\Gamma(SD_{16}, S)$ is isomorphic to the complete bipartite graph $K_{2,8}$. See Fig. 5.2.

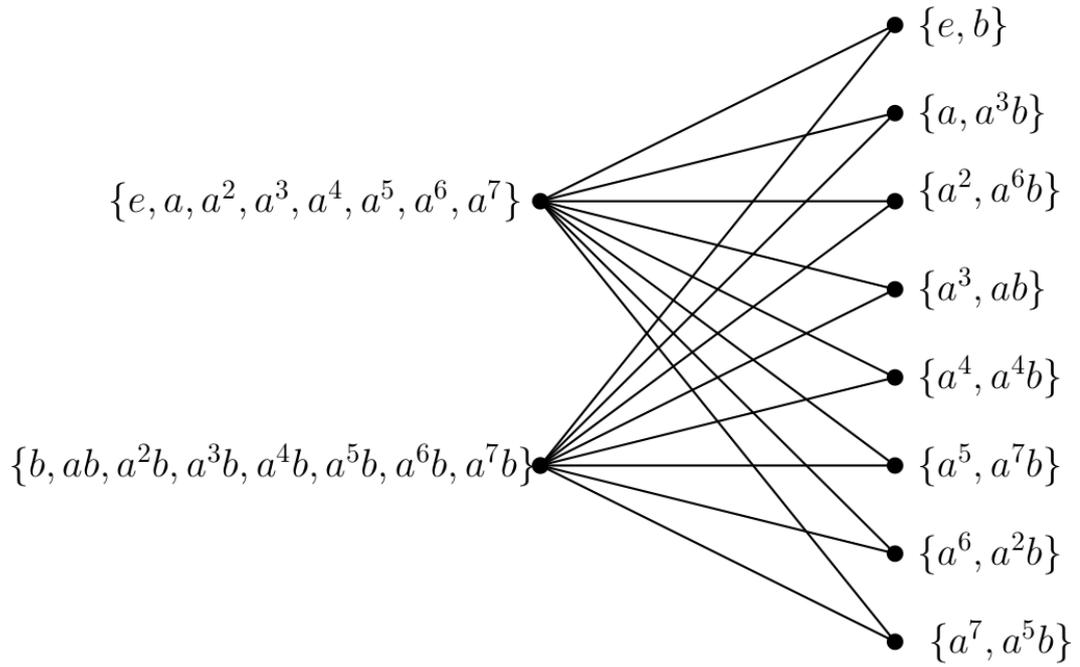

**Figure 5.2:** $\Gamma(SD_{16}, S) \cong K_{2,8}$

# Chapter 6

# Characterisation of finite $G$-graphs

## 6.1 Introduction

Chapter 2 gives results which allow one to construct the $G$-graph of any group given a generating set. Now we give necessary and sufficient conditions for a given graph to be the $G$-graph $\Gamma(G, S)$ of some group $G$ with generating set $S$. By Theorem 3.0.12, $(G, S)$ is determined up to isomorphism. The cardinality of the group and its generating set can be known, in fact, the order of the elements of the generating set is also known and indeed one can give a presentation of the group.

## 6.2 Main Result

**Theorem 6.2.1.** *Let $\Gamma$ be a finite graph with vertex set $V$ and edge set $E$. Then*

> $\Gamma$ *is $k$-partite with $V = V_1 \cup \cdots \cup V_k$, $V_i \cap V_j = \emptyset$ for all $i \neq j$ such that for all $i$, $|V_i| = m_i$, deg $a_i = n_i$, for all $a_i \in V_i$ and $m_i n_i = \dfrac{2|E|}{k}$,*

*if and only if,*

> $\Gamma$ *is isomorphic to the $G$-graph $\Gamma(G, S)$ of a group $G$ with generating set $S$, where $|S| = k, |G| = \dfrac{2|E|}{k(k-1)}$ and for $s_i \in S$, $o(s_i) = \dfrac{n_i}{(k-1)}$, so that $G$ has presentation $\langle s_1, \ldots, s_k | s_i^{o(s_i)} = e \rangle$.*

*Proof.* Suppose

> $\Gamma$ is $k$-partite with $V = V_1 \cup \cdots \cup V_k$, $V_i \cap V_j = \emptyset$ for all $i \neq j$ such that for all $i$, $|V_i| = m_i$, deg $a_i = n_i$, for all $a_i \in V_i$ and $m_i n_i = \dfrac{2|E|}{k}$,

We construct the group $(G, S)$. Since $\Gamma$ is $k$-partite, by proposition 2.4.1, the cardinality of $S$ is $k$. By proposition 2.4.2, $deg(v) = o(s)(|S| - 1)$, so for all $v \in V_i$, $o(s_i) = \frac{deg(v)}{k-1} = \frac{n_i}{k-1}$. Since the cosets corresponding to the vertices in partition $V_i$ of $\Gamma$ partition $G$ and the cardinality of each coset corresponds to the order of a generator, we have $|G| = m_i \times o(s_i) = \frac{m_i n_i}{k-1} = \dfrac{2|E|}{k(k-1)}$. The group $G$ thus has the presentation $\langle s_1, \ldots, s_k | s_i^{o(s_i)} = e \rangle$.

Conversely, suppose

> $\Gamma$ is isomorphic to the $G$-graph $\Gamma(G, S)$ of a group $G$ with generating set $S$, where $|S| = k$, $|G| = \dfrac{2|E|}{k(k-1)}$ and for $s_i \in S$, $o(s_i) = \dfrac{n_i}{(k-1)}$, so that $G$ has presentation $\langle s_1, \ldots, s_k | s_i^{o(s_i)} = e \rangle$.

Then since $|S| = k$, by proposition 2.4.1, the graph $\Gamma$ is $k$-partite with $V = V_1 \cup \cdots \cup V_k$, $V_i \cap V_j = \emptyset$ for all $i \neq j$. From proposition 2.4.2, $deg(v) = o(s)(|S| - 1)$, so $deg(a_i) = \dfrac{n_i}{(k-1)}(k-1) = n_i$, for all $a_i \in V_i$. There are $\frac{|G|}{o(s_i)}$ distinct cosets labelling the vertices in partition $V_i$ corresponding to cosets of $G$ by $\langle s_i \rangle$. So $|V_i| = \frac{|G|}{o(s_i)} = \dfrac{2|E|}{k(k-1)} \div \dfrac{n_i}{(k-1)} = \frac{2|E|}{kn_i}$. If we let $|V_i| = \frac{2|E|}{kn_i} = m_i$ then $m_i n_i = \dfrac{2|E|}{k}$.

□

We are now in position to determine which graphs are $G$-graphs and which are not.

## 6.3 The Turán Graph

The Turán graph $T(n, r)$ named after Pál Turán is a complete r-partite graph. The set of $n$ vertices are partitioned into $r$ different partitions, that is $V = V_1 \cup ... \cup V_r$ such that $|V| = |V_1| + ... + |V_r| = n$. The cardinalities $|V_i|$ of the partitions may differ by 1 and two vertices are adjacent if and only if they belong to different partition sets [3, 36].

The Turán graph $T(n,r)$ has $(n \bmod r)$ partitions of cardinality $\lceil \frac{n}{r} \rceil$ and $r - (n \bmod r)$ partitions of cardinality $\lfloor \frac{n}{r} \rfloor$. Thus $T(n,r) = K_{\lceil \frac{n}{r} \rceil, \lceil \frac{n}{r} \rceil, ..., \lfloor \frac{n}{r} \rfloor, \lfloor \frac{n}{r} \rfloor}$, a complete r-partite graph.

Each vertex $v \in V$ has degree either $deg(v) = n - \lceil \frac{n}{r} \rceil$ or $deg(v) = n - \lfloor \frac{n}{r} \rfloor$. As such it is regular if $n$ is divisible by $r$.

The number of edges $|E|$ is $\lfloor \frac{(r-1)n^2}{2r} \rfloor$ and chromatic number $\chi = r$.

Fig. 6.1 is the Turán graph $T(13, 4)$.

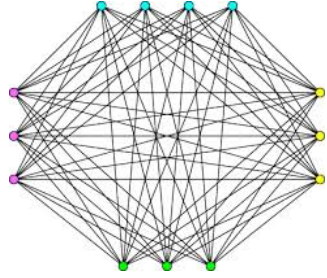

**Figure 6.1:** The Turán graph $T(13, 4)$. Obtained from wikipedia

**Theorem 6.3.1.** *Let $2 \leq r < n$. The Turán graph $T(n,r)$ such that $n$ is divisible by $r$, is a G-graph. [That is, all regular Turán graphs are G-graphs.]*

*Proof.* Since $r$ divides $n$, $\lfloor \frac{n}{r} \rfloor = \lceil \frac{n}{r} \rceil = \frac{n}{r} \in \mathbb{N}$. Thus $T(n,r) = K_{\underbrace{\frac{n}{r}, \frac{n}{r}, ..., \frac{n}{r}}_{r-times}}$. That is for all $i$, $|V_i| = m_i = \frac{n}{r}$ and for all $a_i \in V_i$, $deg(a_i) = n_i = n - \lfloor \frac{n}{r} \rfloor = n - \frac{n}{r} = \frac{n(r-1)}{r}$. Now $m_i n_i = (\frac{n}{r})(\frac{n(r-1)}{r}) = \frac{n^2(r-1)}{r^2}$. Also, $|E| = \frac{(r-1)n^2}{2r}$ so $\frac{2|E|}{k} = \frac{2|E|}{r} = \frac{n^2(r-1)}{r^2}$. Hence the graph satisfies the necessary conditions. The group $G$ has order $|G| = \frac{2|E|}{k(k-1)} = (\frac{n^2(r-1)}{r})(\frac{1}{r(r-1)}) = \frac{n^2}{r^2} = (\frac{n}{r})^2 \in \mathbb{N}$. $o(s_i) = \frac{n_i}{k-1} = (\frac{n(r-1)}{r})(\frac{1}{r-1}) = \frac{n}{r} \in \mathbb{N}$. Thus the group $G$ has $(\frac{n}{r})^2$ elements and $S$ has $r$ elements each with order $\frac{n}{r}$. That is $\langle s_1, s_2, ..., s_r | s_i^{\frac{n}{r}} = e \rangle$. The G-graph $\Gamma(G, S) \cong K_{\underbrace{\frac{n}{r}, \frac{n}{r}, ..., \frac{n}{r}}_{r-times}}$. Hence, the sufficient condition is also satisfied. $\square$

**Example 6.3.2.** The Turán graph $T(6, 3) = K_{2,2,2}$ is the octahedral graph.

**Corollary 6.3.3.** *The Turán graph $T(n,r)$ such that $r^2 = n$ is a G-graph.*

*Proof.* It is straight forward to see that, if $r^2 = n$ then $r$ divides $n$ (specifically $\frac{n}{r} = r$). The graph $T(n,r)$ is thus isomorphic to the G-graph $\Gamma(G, S)$ of the group $G$ with $|G| = n$ and $|S| = r$ such that, for all $s_i \in S$, $o(s_i) = r$. That is $\langle s_1, s_2, ..., s_r | s_i{}^r = e \rangle$. □

**Proposition 6.3.4.** *Let $r > 2$. The Turán graph $T(n,r)$ such that $n$ is not divisible by $r$, is not a G-graph. That is, no non-regular Turán graph $T(n,r)$ $n > 2$ is a G-graph.*

*Proof.* There are vertices with degree $n - \lceil \frac{n}{r} \rceil$ and vertices with degree $n - \lfloor \frac{n}{r} \rfloor$. Since $(n - \lfloor \frac{n}{r} \rfloor) - (n - \lceil \frac{n}{r} \rceil) = \lceil \frac{n}{r} \rceil - \lfloor \frac{n}{r} \rfloor = 1$, for all $a \in V$, $deg(a)$ is odd or even. That is, $V$ has vertices of odd degree and vertices of even degree. On the other hand, $r - 1$ is strictly odd or even and not both. Thus $\frac{deg(a)}{r-1}$ is generally not an integer, but $o(s_i) = \frac{deg(a)}{r-1}$. As such, no such finite group exist. □

**Theorem 6.3.5.** *The Turán graph $T(n,2)$ is a G-graph.*

*Proof.* The case where 2 divides $n$ is satisfied by theorem 6.3.1. In this case, $T(n,2) = K_{\frac{n}{2},\frac{n}{2}}$. Thus the group $G$ has order $|G| = \frac{n^2}{4}$ and $|S| = 2$ such that $o(s_1) = o(s_2) = \frac{n}{2}$. That is $\langle s_1, s_2 | s_1^{\frac{n}{2}} = s_2^{\frac{n}{2}} = e \rangle$.

For the case where 2 does not divide $n$. $T(n,2)$ has one partition of cardinality $\lceil \frac{n}{r} \rceil$ and the other partition of cardinality $\lfloor \frac{n}{r} \rfloor$. Let $\lceil \frac{n}{r} \rceil = x$ and $\lfloor \frac{n}{r} \rfloor = y$, then $T(n,2) = K_{x,y}$, with $x + y = n$. So $|V_1| = m_1 = x$ and for all $a \in V_1$, $deg(a) = n_1 = n - \lceil \frac{n}{r} \rceil = n - x = y$. So $m_1 n_1 = (x)(y) = xy$. But $|E| = xy$ so that $\frac{2|E|}{k} = \frac{2|E|}{r} = \frac{2|E|}{2} = |E| = xy$. Therefore $m_1 n_1 = \frac{2|E|}{k}$. On the other hand, $|V_2| = m_2 = y$ and for all $b \in V_2$, $deg(b) = n_2 = n - \lfloor \frac{n}{r} \rfloor = n - y = x$. So $m_2 n_2 = (y)(x) = xy$. But $|E| = xy$ so that $\frac{2|E|}{k} = \frac{2|E|}{r} = \frac{2|E|}{2} = |E| = xy$. Therefore $m_2 n_2 = \frac{2|E|}{k}$. The necessary condition is thus satisfied.

The group $G$ has order $|G| = \frac{2|E|}{r(r-1)} = \frac{2|E|}{2(1)} = |E| = xy$. The cardinality of $S$ is 2 with $o(s_1) = \frac{n_1}{r-1} = y$ and $o(s_2) = \frac{n_2}{r-1} = x$. That is $|G| = xy$ and $\langle s_1, s_2 | s_1^y = s_2^x = e \rangle$. So the G-graph $\Gamma(G,S)$ is isomorphic to $T(n,2) = K_{x,y}$. □

**Remark 6.3.6.** $T(n,2)$ is a complete bipartite graph.

**Example 6.3.7.** The G-graph $\Gamma(\mathbb{Z}/x\mathbb{Z} \times \mathbb{Z}/y\mathbb{Z}, S)$ where $S = \{(1,0),(0,1)\}$ is isomorphic to $K_{x,y}$.

**Corollary 6.3.8.** *The Moore graph is a G-graph.*

**Remark 6.3.9.** The Moore graph is the Turán graph $T(n,2)$ with $n$ even.

## 6.4 The Turán graph $T(n, n)$ and the Trivial group

We note that, the Turán graph $T(n,n)$ is the complete regular graph $K_n$. On the other hand, the trivial group has only the identity as an element. The following result also appears in [34].

**Proposition 6.4.1.** *Let $G$ be the trivial group, i.e., $G = \{e\}$ and let $S$ be an $n$-multi set of the identity element, (i.e., $S = \{e, e, e, ..., e\}$ such that $|S| = n$). Then the $G$-graph $\Gamma(G, S) \cong T(n, n) = K_n$.*

*Proof.* Since $|S| = n$, the $G$-graph is $n$-partite. Moreover, for all $s_i \in S$, $o(s) = 1$ so the coset of each cyclic subgroup $\langle e \rangle$ is $\{e\}$. Thus each of the $n$ partition sets has just one vertex, $\{e\}$, intersecting the other $n-1$ vertices (cosets) $\{e\}$. Hence we have the complete regular graph. $\square$

**Example 6.4.2.** If $G = \{e\}$ and $S = \{e, e, e, e\}$, then $\Gamma(G, S) \cong T(4, 4) = K_4$, the tetrahedral graph.

## 6.5 The Platonic graphs

The platonic graphs are the skeletons of the platonic solids. There are five platonic graphs. We have seen that, the octahedral and the cubical graphs are $G$-graphs. Moreover, we have seen that, the tetrahedral graph is a $G$-graph. We show below that they are the only platonic graphs that are $G$-graphs. The other platonic graphs are the icosahedral graph and the dodecahedral graph.

### 6.5.1 The Icosahedral graph

The icosahedral graph has $|V| = 12$ and $|E| = 30$. It is regular of degree 5 and has chromatic number $\chi = 4 = k$. Though this graph satisfies the necessary condition, it fails the sufficient condition. That is, there exists a partition $V = V_1 \cup V_2 \cup V_3 \cup V_4$ such that, for all $V_i$, $|V_i| = m_i = 3$ and for all $a_i \in V_i$ $deg(a_i) = n_i = 5$ so that, $m_i n_i = 15 = \frac{2|E|}{k}$. However, the corresponding group of order $\frac{2|E|}{k(k-1)} = \frac{60}{12} = 5$ must be generated by $S$ such that $|S| = 4$ and each $s_i \in S$ has order $\frac{n_i}{k-1} = \frac{5}{3}$ which is impossible.

### 6.5.2 The Dodecahedral graph

The dodecahedral graph has $|V| = 20$ and $|E| = 30$. It is regular of degree 3 and has chromatic number $\chi = 3 = k$. It is straight forward to see that, we can not get $(G, S)$, a group $G$ of order $\frac{2|E|}{k(k-1)} = \frac{60}{6} = 10$ generated by $S$ such that $|S| = 3$ and each $s_i \in S$ has order $\frac{n_i}{k-1} = \frac{3}{2}$.

## 6.6 Biregular and Regular Bipartite Graphs

A semi-regular bipartite graph also known as a biregular graph is a bipartite graph such that vertices in the same bipartition have the same degree. Thus the bipartite graph $H = (U \cup V, E)$ is $(m, n)$-biregular if for all $u_i \in U$, $deg(u_i) = m$ and for all $v_i \in V$, $deg(v_i) = n$ [1, 20, 32]. Therefore every regular bipartite graph is as well biregular.

**Remark 6.6.1.** [1] An $(m, n)$-biregular graph $H = (U \cup V, E)$ satisfies the equation:

$$m|U| = n|V| = |E|.$$

The following theorem gives a large class of graphs that are $G$-graphs.

**Theorem 6.6.2.** *All $(m, n)$-biregular graphs are $G$-graphs.*

*Proof.* Let $H = (U \cup V, E)$ be an $(m, n)$-biregular graph. Then by remark 6.6.1, $m|U| = n|V| = |E|$ satisfying the necessary conditions of theorem 6.2.1. Now since $H$ is bipartite, $|S| = k = 2$. We can get a group $G$ such that $|G| = \frac{2|E|}{k(k-1)} = \frac{2|E|}{2} = |E|$. Moreover for $s_1 \in S$, $o(s_1) = \frac{m}{2-1} = m$ and $s_2 \in S$ $o(s_2) = \frac{n}{2-1} = n$. Thus an $(m, n)$-biregular graph $H = (U \cup V, E)$ is isomorphic to the $G$-graph $\Gamma(G, S)$, where the couple $(G, S)$ is such that $|G| = |E|$ and $S = \{s_1, s_2\}$ with $o(s_1) = m$, $o(s_2) = n$. That is the group of order $|E|$ with presentation $\langle s_1, s_2 | s_1^m = s_2^n = e \rangle$. □

**Example 6.6.3.** The star graph $S_n$ is a $(1, n)$-biregular graph. It is isomorphic to $\Gamma(G, S)$ where $|G| = n$ and $S = \{s_1, s_2\}$ with the presentation $\langle s_1, s_2 | s_1^1 = s_2^n = e \rangle = \langle e, s_2 | s_2^n = e \rangle$. In particular, $G$ is the cyclic group of order n, $C_n = \{e, a, a^2, a^3, a^4, ..., a^n\}$ and $S = \{e, a\}$.

**Example 6.6.4.** The Rhombic Dodecahedral graph, a $(3,4)$-biregular graph is a $G$-graph. In fact it is isomorphic to $\Gamma(G, S)$ where $|G| = 24$ and $S = \{s_1, s_2\}$ with the presentation $\langle s_1, s_2 | s_1^3 = s_2^4 = e \rangle$.

Since every regular bipartite graph is also biregular, corollary 6.6.5 is immediate.

**Corollary 6.6.5.** *Every $n$-regular bipartite graph is a $G$-graph. In which case the group $G$ is of order $|E|$ with presentation $\langle s_1, s_2 | s_1^n = s_2^n = e \rangle$.*

A graph $\Gamma$ is said to be *vertex-transitive*, if for any pair of vertices $v, u \in V$, there exists an automorphism $\sigma \in Aut(\Gamma)$ such that $\sigma(v) = u$.

Similarly, a graph $\Gamma$ is said to be *edge-transitive*, if for any pair of edges $e, e' \in E$, there exists an automorphism $\sigma \in Aut(\Gamma)$ such that $\sigma(e) = e'$.

**Remark 6.6.6.** [1, 27] Except for graphs with isolated vertices, every edge-transitive graph that is not also vertex-transitive is biregular.

As a result, all such classes of graphs are $G$-graphs.

**Remark 6.6.7.** [2, 16] The Levi graphs of geometric configurations are biregular.

As such the following Levi graphs are $G$-graphs. One can also refer to [13] for their corresponding groups.

- The Desargues graph
- The Heawood graph
- The Möbius–Kantor graph
- The Pappus graph
- The Gray graph
- The Ljubljana graph
- The Tutte eight-cage. $|G| = 45$ and $\langle s_1, s_2 | s_1^3 = s_2^3 = e \rangle$.
- The four-dimensional hypercube graph $Q_4$.

The Tutte eight-cage also called the Tutte-Coxeter graph is a Moore graph shown to be $G$-graph in Corollary 6.3.8.

## 6.7 The n-dimensional Hypercube Graph

The $n$-cube, denoted by $Q_n$ is a graph representing the connections between vertices in an $n$-dimensional cube. It is bipartite [23]. It is $n$-regular, has $2^n$ vertices and $2^{n-1}n$ edges.

**Proposition 6.7.1.** *The $n$-cube is a $G$-graph.*

*Proof.* Let $H = (V_r \cup V_b, E)$ be an $n$-dimensional cube graph, then $|V_r| = |V_b| = 2^{n-1}$. It is regular of degree $n$ so the necessary conditions of theorem 6.2.1 are satisfied since $n.2^{n-1} = |E|$. Moreover, the sufficient conditions are satisfied since we can get a group $G$ such that $|G| = 2^{n-1}n$ and $S = \{s_1, s_2\}$ such that for all $s_i \in S$, $o(s_i) = \frac{n}{2-1} = n$. That is, there exists a group $G$ of order $2^{n-1}n$ with the presentation $\langle s_1, s_2 | s_1^n = s_2^n = e \rangle$. $\square$

**Remark 6.7.2.** Since the $n$-cube is bipartite and regular, proposition 6.7.1 is satisfied in corollary 6.6.5.

**Example 6.7.3.**

- $Q_3$ is isomorphic to $\Gamma(G, S)$ where $|G| = 12$ and $s = \{s_1, s_2\}$ with the presentation $\langle s_1, s_2 | s_1^3 = s_2^3 = e \rangle$.

- $Q_4$ is isomorphic to $\Gamma(G, S)$ where $|G| = 32$ and $s = \{s_1, s_2\}$ with the presentation $\langle s_1, s_2 | s_1^4 = s_2^4 = e \rangle$.

- $Q_5$ is isomorphic to $\Gamma(G, S)$ where $|G| = 80$ and $s = \{s_1, s_2\}$ with the presentation $\langle s_1, s_2 | s_1^5 = s_2^5 = e \rangle$.

- $Q_6$ is isomorphic to $\Gamma(G, S)$ where $|G| = 192$ and $s = \{s_1, s_2\}$ with the presentation $\langle s_1, s_2 | s_1^6 = s_2^6 = e \rangle$.

## 6.8 Characterisation of Bipartite $G$-graphs

The following theorem is in agreement with theorem 6.1 in [11].

**Theorem 6.8.1.** *A connected bipartite graph $\Gamma$ is a $G$-graph if and only if it is biregular.*

*Proof.* Suppose $\Gamma$ is a bipartite $G$-graph. Then there exists the couple $(G, S)$ such that by propostion 2.4.1, $|S| = 2$ and by proposition 2.4.5, vertices in the same bipartition have the same degree. Thus $\Gamma$ is biregular.

Conversely, suppose $\Gamma$ is biregular, then by theorem 6.6.2, $\Gamma$ is a $G$-graph. $\square$

**Remark 6.8.2.** In the above theorem, the elements in $S$ need not be unique, except that $|S| = 2$. We also note that, we are dealing with loopless graphs.

## Chapter 7

# The Spectrum of a finite $G$-graph and some infinite $G$-graphs

In this chapter, we consider some equally important aspects of $G$-graphs which will require further investigation. We present some preliminary results. We consider the adjacency matrix and the spectrum of a finite $G$-graph. Then we shift our focus to infinite $G$-graphs.

## 7.1 Adjacency Matrix of a finite $G$-graph

Let $\Gamma(V, E)$ be a graph where $V = \{1, 2, 3, ..., n\}$. For $1 \leq i, j \leq n$, the adjacency matrix $M(\Gamma) = (a_{i,j})$ is defined as

$$(a_{i,j}) = \begin{cases} k & \text{if } (i,j) \in E \\ 0 & \text{if } (i,j) \notin E \end{cases}$$

where $k$ is the number of edges between $i$ and $j$.

We extend this definition to $G$-graphs.

**Definition 7.1.1.** The adjacency matrix $M(\Gamma)$ of the $G$-graph $\Gamma = \Gamma(G, S)$ is the matrix with rows and columns indexed by the distinct cosets of the cyclic subgroups $\langle s_i \rangle$, $s_i \in S$ such that $(a_{i,j}) = |\langle s_i \rangle x \cap \langle s_j \rangle y|$ for all $i \neq j$ and $(a_{i,i}) = 0$, where $x, y \in G$.

Results in previous chapters enable us to make the following observation about the adjacency matrix of a finite $G$-graph. Now since, the $G$-graph $\Gamma = \Gamma(G, S)$ of the group $G$ with respect to the generating set $S$ is $|S| = k-$partite, there exists a permutation matrix $P$ such that, the adjacency matrix $M(\Gamma)$ is given by

$$M(\Gamma) = \begin{pmatrix} \mathbf{0} & X_1 & X_2 & \ldots & X_k \\ X_1^T & \mathbf{0} & Y_2 & \ldots & Y_k \\ X_2^T & Y_2^T & \mathbf{0} & \ldots & Z_k \\ \vdots & \vdots & \vdots & \ldots & \vdots \\ X_k^T & Y_k^T & Z_k^T & \ldots & \mathbf{0} \end{pmatrix}$$

where $\mathbf{0}, X_i, Y_i, Z_i$ etc, are block matrices. The $\mathbf{0}$'s are $M_{n_i \times n_i}$ null matrices along the diagonal. Thus $|M_{n_i \times n_i}| = \chi(\Gamma)$ and $n_i = |V_i|$, the number of vertices in that partition.

For a bipartite $G$-graph $\Gamma$, we have

$$\begin{pmatrix} \mathbf{0} & X \\ X^T & \mathbf{0} \end{pmatrix}$$

And for a tri-partite $G$-graph $\Gamma$, we have

$$\begin{pmatrix} \mathbf{0} & X_1 & X_2 \\ X_1^T & \mathbf{0} & Y_2 \\ X_2^T & Y_2^T & \mathbf{0} \end{pmatrix}$$

**Example 7.1.2.** The adjacency matrix of $\Gamma(D_{10}, \{s, t\})$ is

$$\begin{pmatrix} 0 & 0 & 1 & 1 & 1 & 1 & 1 \\ 0 & 0 & 1 & 1 & 1 & 1 & 1 \\ 1 & 1 & 0 & 0 & 0 & 0 & 0 \\ 1 & 1 & 0 & 0 & 0 & 0 & 0 \\ 1 & 1 & 0 & 0 & 0 & 0 & 0 \\ 1 & 1 & 0 & 0 & 0 & 0 & 0 \\ 1 & 1 & 0 & 0 & 0 & 0 & 0 \end{pmatrix}$$

**Example 7.1.3.** The adjacency matrix of $\Gamma(S_3, \{(12), (13), (23)\})$ is

$$\begin{pmatrix} 0 & 0 & 0 & 1 & 1 & 0 & 1 & 1 & 0 \\ 0 & 0 & 0 & 1 & 0 & 1 & 0 & 1 & 1 \\ 0 & 0 & 0 & 0 & 1 & 1 & 1 & 0 & 1 \\ 1 & 1 & 0 & 0 & 0 & 0 & 1 & 0 & 1 \\ 1 & 0 & 1 & 0 & 0 & 0 & 0 & 1 & 1 \\ 0 & 1 & 1 & 0 & 0 & 0 & 1 & 1 & 0 \\ 1 & 0 & 1 & 1 & 0 & 1 & 0 & 0 & 0 \\ 1 & 1 & 0 & 0 & 1 & 1 & 0 & 0 & 0 \\ 0 & 1 & 1 & 1 & 1 & 0 & 0 & 0 & 0 \end{pmatrix}$$

### 7.1.4 Properties

The adjacency matrix of a $G$-graph has the general properties of an adjacency matrix of any graph. We make the following observations in connection to the groups they represent. We show how to identify when a graph is not a $G$-graph from its adjacency matrix.

- The sum $d_s(n)$ of the entries of a column/row gives the degree of the corresponding vertex. If $|M_{n_i \times n_i}| = 2$, then $d_s(n) = o(s)$, the order of the corresponding generator. If $|M_{n_i \times n_i}| = l > 2$, then the order of the corresponding generator is $o(s) = \frac{d_s(n)}{l-1}$. Moreover, the order of the group $|G| = (n_i)o(s)$.

- All columns/rows of $M(\Gamma)$ in the block $M_{n_i \times n_i}$ have the same $d_s(n)$. This is due to Proposition 2.4.5.

- Let $|M_{n_i \times n_i}| > 2$. If there exist columns/rows $v_1, v_2$ such that $d_s(v_1) - d_s(v_2) = 1$, then the graph represented by the adjacency matrix $M(\Gamma)$ is not a $G$-graph.

- The number of edges $|E| = \frac{\sum d_s(n)}{2}$, the total sum of the entries of the upper or lower diagonal. This is because, $M(\Gamma)$ is symmetric.

## 7.2 The Spectrum of a $G$-graph

The set of eigenvalues of the adjacency matrix of a graph is called the spectrum of the graph. The spectrum of a graph is an important graph invariant. It has several applications. See [14, 33]. One of its application is the energy of a graph. According to Majstorović [29], there is interest in number of distinct eigenvalues and in the energy of a graph. The *energy* $E(G)$ of a simple graph is the sum of the absolute values of its eigenvalues. The definition of graph energy is closely connected to the $\pi$-electron energy of a conjugated carbon molecule, computed using the H$\ddot{u}$ckel theory. So results on graph energy assume special significance. A graph $G$ with $n$ vertices is said to be "hyperenergetic" if its energy $E(G) > 2n - 2$, and to be "hypoenergetic" if its energy $E(G) < n$. Check [5, 25, 28] for details. In this section we compute the spectrum of the adjacency matrix and energy for some $G$-graphs. Exponents imply multiplicities. We used *Sage* for our computations.

1. Group $G$ is the Klein group. $G = \{e, a, b, ab\}$. The $G$-graph $\Gamma(G, \{a, b, ab\}) \cong$ the Octahedral graph.

   Spectrum $=\{0^3, -2^2, 4\}$. It has 3 distinct eigenvalues. The energy $E(G) \geq n$, the number of vertices.

2. Group $G$ is the Dihedral group. $G = D_{10}$. The $G$-graph $\Gamma(G, \{a, b\}) \cong K_{2,5}$, where $o(a) = 2$ and $o(b) = 5$.

   Spectrum $=\{0^5, -3.16228, 3.16228\}$. It has 3 distinct eigenvalues. The energy $E(G) < n$, the number of vertices.

3. Group $G$ is the Dihedral group. $G = D_4$. The $G$-graph $\Gamma(G, \{a, b\}) \cong C_4$, where $o(a) = 2$ and $o(b) = 2$.

   Spectrum $=\{0^2, -1, 1\}$. It has 3 distinct eigenvalues. The energy $E(G) < n$, the number of vertices.

4. Group $G$ is the generalized quaternion group $Q_3$. The $G$-graph $\Gamma(G, \{a, b\}) \cong K_{2,3}^2$, where $o(a) = 4$ and $o(b) = 6$.

   Spectrum $=\{0^3, -4.89898, 4.89898\}$. It has 3 distinct eigenvalues.

5. Group $G$ is the quaternion group. $G = Q$. The $G$-graph $\Gamma(G, \{a, b\}) \cong K_{2,2}^2$, where $o(a) = 4$ and $o(b) = 4$.

   Spectrum $=\{0^3, -4, 4\}$. It has 3 distinct eigenvalues.

6. Group $G$ is the Abelian group $G = C_3 \times C_3$. The $G$-graph $\Gamma(G, \{(1, 0), (0, 1)\}) \cong K_{3,3}$, where $o(a) = 3$ and $o(b) = 3$.

   Spectrum $=\{0^4, -3, 3\}$. It has 3 distinct eigenvalues. The energy $E(G) = n$, the number of vertices.

Generally, for $G = C_k \times C_k$. The $G$-graph $\Gamma(G, \{(1,0), (0,1)\}) \cong K_{k,k}$, where $o(a) = k$ and $o(b) = k$.

Spectrum $=\{0^{2k-2}, -k, k\}$. It has 3 distinct eigenvalues. Its characteristic polynomial is $(\lambda - k)(\lambda + k)(\lambda)^{2k-2}$. The energy $E(G) = n$, the number of vertices.

7. Group $G$ is the Trivial group $G = \{e\}$. The $G$-graph $\Gamma(G, \{e, e, e, e, e\}) \cong K_5 =$Turan(5,5).

    Spectrum $=\{-1^4, 4\}$. It has 2 distinct eigenvalues. The energy $E(G) > n$, the number of vertices.

    Generally, for $K_n$, spectrum is $\{-1^{n-1}, (n-1)\}$. The energy $E(G) = 2n-2 > n$, the number of vertices.

8. Group $G$ is the symmetric group. $G = S_3$. The $G$-graph $\Gamma(G, \{(12), (13), (23)\})$.

    Spectrum $=\{-2^4, 1^4, 4\}$. It has 3 distinct eigenvalues. The energy $E(G) > n$, the number of vertices.

9. Group $G$ is the symmetric group. $G = S_4$. The $G$-graph $\Gamma(G, \{(12), (234)\})$.

    Spectrum $=\{-2.449, -2^3, -1.4142^3, 1.4142^3, 2^3, 0^6, 2.449\}$. It has 7 distinct eigenvalues. The energy $E(G) > n$, the number of vertices.

10. Group $G$ is the alternating group. $G = A_4$. The $G$-graph $\Gamma(G, \{(123), (134)\})$.

    Spectrum $=\{-3, -1^3, 1^3, 3\}$. It has 4 distinct eigenvalues. The energy $E(G) > n$, the number of vertices.

It is known that the energy $E(G) \geq n$, the number of vertices, is satisfied by

i. regular graphs

ii. graphs with no zero eigenvalues [as in 7 and 8 above] and

iii. graphs satisfying edges $\geq \frac{n^2}{4}$.

We were unable to discern any particular properties from the spectrum. In most cases, the energy $E(G) \geq n$, the number of vertices. We only found hypoenergetic $G$-graphs in 2 and 3.

## 7.3 Infinite $G$-graphs

We study finitely generated infinite groups. An infinite graph has an infinite number of vertices or edges or both. We note that, a graph is locally finite if each of its vertices has finite degree. We consider the special linear group $SL_2(\mathbb{Z})$ and an infinite non-Abelian matrix group.

### 7.3.1 $G$-graph of $SL_2(\mathbb{Z})$

The group $SL_2(\mathbb{Z})$ consists of $2 \times 2$ integer matrices with determinant 1 under multiplication. The group is infinite, and is generated by $s_1 = \begin{pmatrix} 0 & -1 \\ 1 & 0 \end{pmatrix}$ and $s_2 = \begin{pmatrix} 1 & 1 \\ 0 & 1 \end{pmatrix}$ with $o(s_1) = 4$ and $o(s_2) = \infty$ since $s_1^k = \begin{pmatrix} 1 & k \\ 0 & 1 \end{pmatrix}$. We also note that, $x = s_1 s_2 = s_1 = \begin{pmatrix} 0 & -1 \\ 1 & 1 \end{pmatrix}$ with $o(x) = 6$. Thus $SL_2(\mathbb{Z})$ is an infinite group that can be generated by two elements of finite order [17].

**The $G$-graph $\Gamma(SL_2(\mathbb{Z}), \{s_1, x\})$**

We shall explicitly construct the graph.

We find the cosets of the cyclic subgroup $\langle s_1 \rangle$.

We note the coset that contains the identity element is

$$\langle s_1 \rangle e = \left\{ \begin{pmatrix} 1 & 0 \\ 0 & 1 \end{pmatrix}, \begin{pmatrix} 0 & -1 \\ 1 & 0 \end{pmatrix}, \begin{pmatrix} -1 & 0 \\ 0 & -1 \end{pmatrix}, \begin{pmatrix} 0 & 1 \\ -1 & 0 \end{pmatrix} \right\}.$$

Now for all $\begin{pmatrix} a & b \\ c & d \end{pmatrix} \in SL_2(\mathbb{Z})$, the coset of the cyclic subgroup $\langle s_1 \rangle$ is of the form

$$\langle s_1 \rangle \begin{pmatrix} a & b \\ c & d \end{pmatrix} = \left\{ \begin{pmatrix} a & b \\ c & d \end{pmatrix}, \begin{pmatrix} -c & -d \\ a & b \end{pmatrix}, \begin{pmatrix} -a & -b \\ -c & -d \end{pmatrix}, \begin{pmatrix} c & d \\ -a & -b \end{pmatrix} \right\}.$$

Thus there are infinitely many cosets $\langle s_1 \rangle g$, $g \in SL_2(\mathbb{Z})$. The cardinality of each coset is 4.

Now we find the cosets of the cyclic subgroup $\langle x \rangle$.

We note the coset that contains the identity element is

$$\langle x \rangle e = \left\{ \begin{pmatrix} 1 & 0 \\ 0 & 1 \end{pmatrix}, \begin{pmatrix} 0 & -1 \\ 1 & 1 \end{pmatrix}, \begin{pmatrix} -1 & -1 \\ 1 & 0 \end{pmatrix}, \begin{pmatrix} -1 & 0 \\ 0 & -1 \end{pmatrix}, \begin{pmatrix} 0 & 1 \\ -1 & -1 \end{pmatrix}, \begin{pmatrix} 1 & 1 \\ -1 & 0 \end{pmatrix} \right\}.$$

Now for all $\begin{pmatrix} a & b \\ c & d \end{pmatrix} \in SL_2(\mathbb{Z})$, the coset of the cyclic subgroup $\langle x \rangle$ is of the form

$$\langle x \rangle \begin{pmatrix} a & b \\ c & d \end{pmatrix} = \left\{ \begin{pmatrix} a & b \\ c & d \end{pmatrix}, \begin{pmatrix} -c & -d \\ a+c & b+d \end{pmatrix}, \begin{pmatrix} -a-c & -b-d \\ a & b \end{pmatrix}, \right.$$

$$\left. \begin{pmatrix} -a & -b \\ -c & -d \end{pmatrix}, \begin{pmatrix} c & d \\ -a-c & -b-d \end{pmatrix}, \begin{pmatrix} a+c & b+d \\ -a & -b \end{pmatrix} \right\}.$$

Thus there are infinitely many cosets $\langle x \rangle g$, $g \in SL_2(\mathbb{Z})$. The cardinality of each coset is 6.

We note that, for all cosets $\langle s_1 \rangle u$, $u \in SL_2(\mathbb{Z})$ and $\langle x \rangle v$, $v \in SL_2(\mathbb{Z})$, $|\langle s_1 \rangle u \cap \langle x \rangle v| = 0$ or 2. Thus the graph is 2-edged and hence every vertex in the partition corresponding to $\langle s_1 \rangle u$, $u \in SL_2(\mathbb{Z})$ is adjacent to two vertices in the partition corresponding to $\langle x \rangle v$, $v \in SL_2(\mathbb{Z})$. So every vertex in the partition corresponding to $\langle x \rangle v$, $v \in SL_2(\mathbb{Z})$ is adjacent to three vertices in the partition corresponding to $\langle s_1 \rangle u$, $u \in SL_2(\mathbb{Z})$. We conclude with the following;

**Proposition 7.3.2.** *The $G$-graph $\Gamma(SL_2(\mathbb{Z}), \{s_1, x\})$ is isomorphic to the infinite double edged bipartite graph $(U \cup V, E)$ such that $|U| = \infty = |V|$, for all $u_i \in U$, $deg(u_i) = 4$, for all $v_i \in V$, $deg(v_i) = 6$ and $|E| = \infty$.*

Fig. 7.1 is the $G$-graph $\Gamma(SL_2(\mathbb{Z}), \{s_1, x\})$. Vertices in the first partition have degree 4 while vertices in the second partition have degree 6.

### 7.3.3 The $G$-graph of an infinite non-Abelian Matrix Group

The infinite non-Abelian matrix group $G = \left\{ \begin{pmatrix} a & b \\ 0 & 1 \end{pmatrix} : a = \pm 1, b \in \mathbb{Z} \right\}$ is finitely generated. The group has generating set $S = \left\{ \begin{pmatrix} -1 & 0 \\ 0 & 1 \end{pmatrix}, \begin{pmatrix} 1 & 1 \\ 0 & 1 \end{pmatrix} \right\}$.

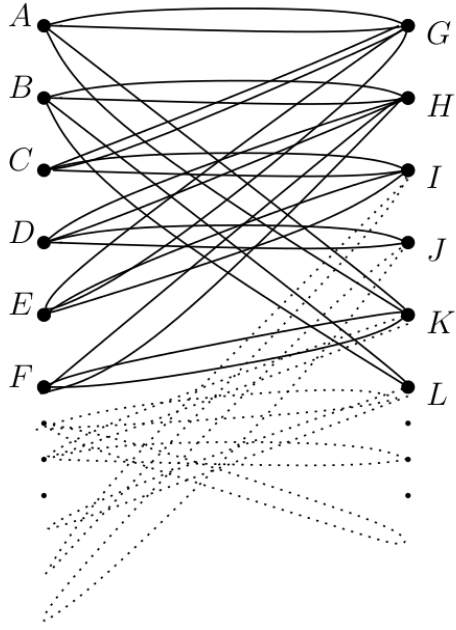

**Figure 7.1:** $\Gamma(G, \{s_1, x\})$

Let $s_0 = \begin{pmatrix} -1 & 0 \\ 0 & 1 \end{pmatrix}$ and $s_1 = \begin{pmatrix} 1 & 1 \\ 0 & 1 \end{pmatrix}$. Then $o(s_0) = 2$ and $o(s_1) = \infty$ [17]. We also note that, $s_2 = s_0 s_1 = \begin{pmatrix} -1 & 0 \\ 0 & 1 \end{pmatrix} \begin{pmatrix} 1 & 1 \\ 0 & 1 \end{pmatrix} = \begin{pmatrix} -1 & -1 \\ 0 & 1 \end{pmatrix}$ and $o(s_2) = 2$. Thus since $G = \langle s_0, s_1 \rangle$, we can also say $G = \langle s_0, s_2 \rangle$.

**The $G$-graph $\Gamma(G, \{s_0, s_2\})$**

We shall explicitly construct the graph.

We find the cosets of the cyclic subgroup $\langle s_0 \rangle$.

We note the coset that contains the identity element is

$$\langle s_0 \rangle e = \left\{ \begin{pmatrix} 1 & 0 \\ 0 & 1 \end{pmatrix}, \begin{pmatrix} -1 & 0 \\ 0 & 1 \end{pmatrix} \right\}.$$

Now for all $\begin{pmatrix} a & b \\ 0 & 1 \end{pmatrix} \in G$, the coset of the cyclic subgroup $\langle s_0 \rangle$ is of the form

$$\langle s_0 \rangle \begin{pmatrix} a & b \\ 0 & 1 \end{pmatrix} = \left\{ \begin{pmatrix} a & b \\ 0 & 1 \end{pmatrix}, \begin{pmatrix} -a & -b \\ 0 & 1 \end{pmatrix} \right\}.$$

Thus there are infinitely many cosets $\langle s_0 \rangle g$, $g \in G$. The cardinality of each coset is 2.

Now we find the cosets of the cyclic subgroup $\langle s_2 \rangle$.

We note the coset that contains the identity element is

$$\langle s_2 \rangle e = \left\{ \begin{pmatrix} 1 & 0 \\ 0 & 1 \end{pmatrix}, \begin{pmatrix} -1 & -1 \\ 0 & 1 \end{pmatrix} \right\}.$$

Now for all $\begin{pmatrix} a & b \\ 0 & 1 \end{pmatrix} \in G$, the coset of the cyclic subgroup $\langle s_2 \rangle$ is of the form

$$\langle s_2 \rangle \begin{pmatrix} a & b \\ 0 & 1 \end{pmatrix} = \left\{ \begin{pmatrix} a & b \\ 0 & 1 \end{pmatrix}, \begin{pmatrix} -a & -b-1 \\ 0 & 1 \end{pmatrix} \right\}.$$

Thus there are infinitely many cosets $\langle s_2 \rangle g$, $g \in G$. The cardinality of each coset is 2.

We note that, for all cosets $\langle s_0 \rangle u$, $u \in G$ and $\langle s_2 \rangle v$, $v \in G$ $|\langle s_0 \rangle u \cap \langle s_2 \rangle v| = 0$ or $1$. Thus the graph is single edged and hence every vertex in the partition corresponding to $\langle s_0 \rangle u$, $u \in G$ is adjacent to one vertex in the partition corresponding to $\langle s_2 \rangle v$, $v \in G$. So every vertex in the partition corresponding to $\langle s_2 \rangle v$, $v \in G$ is adjacent to one vertices in the partition corresponding to $\langle s_0 \rangle u$, $u \in G$. We conclude with the following;

**Proposition 7.3.4.** *The $G$-graph $\Gamma(G, \{s_0, s_2\})$ is isomorphic to the infinite 2-regular bipartite graph $(U \cup V, E)$ such that $|U| = \infty = |V|$ and $|E| = \infty$. In other words, the $G$-graph $\Gamma(G, \{s_0, s_2\})$ is isomorphic to the 2-way infinite Path graph.*

For a general element $x \in G$, say $x = \begin{pmatrix} a & b \\ 0 & 1 \end{pmatrix}$, Fig 7.2 is the $G$-graph $\Gamma(G, \{s_0, s_2\})$.

The Path graph $P_\infty$, is $...A, B, C, D, E, F, G, H, ....$ See in Fig. 7.3.

**Remark 7.3.5.** The Cayley graph of an infinite group with respect to a finite generating set is always locally finite, but the $G$-graph of an infinite group with respect to a finite generating set may or may not be locally finite. We have

**Proposition 7.3.6.** *Let $G$ be an infinite group and suppose $S = \{s_1, s_2, ..., s_k\}$ generates $G$. The $G$-graph $\Gamma(G, S)$ is locally finite if and only if $o(s_i) \neq \infty$ for all $s_i \in S$.*

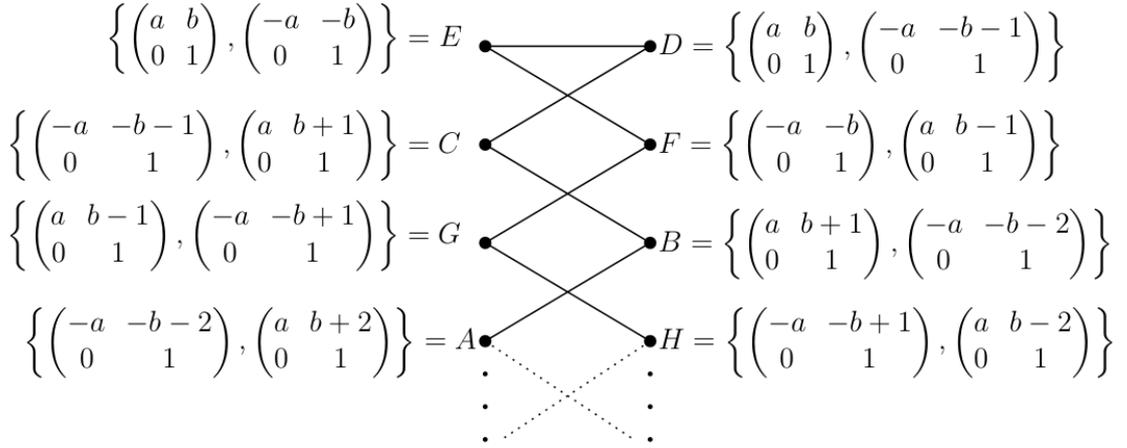

Figure 7.2: $\Gamma(G, \{s_0, s_2\}) \cong P_\infty$

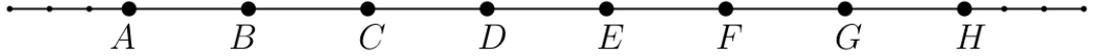

$\quad\quad A\quad\quad B\quad\quad C\quad\quad D\quad\quad E\quad\quad F\quad\quad G\quad\quad H$

Figure 7.3: $\Gamma(G, \{s_0, s_2\}) \cong P_\infty$

*Proof.* Suppose $o(s_i) \neq \infty$ for all $s_i \in S$. Then the proof is from the fact that, the vertices of the graph are the cosets of the generators. Since the order of a generator is finite, the corresponding coset is also finite. The degree of each vertex is thus $o(s_i)(|S| - 1)$. Hence $\Gamma(G, S)$ is locally finite.

On the other hand, if $\Gamma(G, S)$ is locally finite then the degree $deg(v_{s_i}) = o(s_i)(|S|-1)$ of each vertex $\langle s_i \rangle x$ is finite. Hence $o(s_i) \neq \infty$ for all $s_i \in S$. □

If we suppose $o(s_i) = \infty$ for some $s_i \in S$. The cosets $\langle s_i \rangle g$, $g \in G$ of $s_i$ will be of infinite order. Thus the vertices corresponding to each of these cosets will be of infinite degree. Hence if there exists an $s_i \in S$ such that $o(s_i) = \infty$, then the G-graph $\Gamma(G, S)$ is not locally finite.

# Chapter 8

# Conclusion and Further Works

Cayley graphs have important structural properties that make them very useful. However, some of their properties (for instance, that Cayley graphs are always regular) pose limitations on their application. This prompted the introduction of $G$-graphs as an alternate way of representing groups as graphs. Cayley graphs have been extensively studied. The study of $G$-graphs on the other hand is still a young and rich area to be explored. In this thesis, we studied the structure of $G$-graphs and presented results that will help their identification and characterisation. There are other important characterisations already presented by others which we have not addressed. These results use the automorphisms on the graph. Our results are mainly based on the structure of the graph and its connection to the group. In fact, our results emphasise the fact that, one can easily get the $G$-graph of a group by merely knowing the order of the group, the cardinality of the generating set and the orders of the generators.

## 8.1 Main Contributions

Our main contributions to the area of $G$-graphs include the following:

1. We made some useful contributions toward the properties of $G$-graphs. See chapter 2.

2. We presented results that show *how* the $G$-graph of a group depends on the generating sets of the group. See chapter 3.

3. For various generating sets for both the Symmetric group $S_n$ and the alternating group $A_n$, we gave the corresponding $G$-graphs. We also gave the $G$-graphs of the quaternion, the dihedral group generated by involutions and the semi-dihedral group. See chapter 4 and 5.

4. One major contribution was the presentation of a characterisation for finite $G$-graphs. Unlike previous characterisations, our result emphasised the structure of the graph and its connection to the group. In fact, we gave necessary and sufficient conditions for a graph to be a $G$-graph. This characterisation allows for repetition of elements of the generating set thereby capturing wider class of graphs. Using this result, we showed certain family of graphs that are $G$-graphs and those which are not. This approach can be applied to all graphs. We also presented a characterisation of bipartite $G$-graphs, which states that a connected bipartite graph $\Gamma$ is a $G$-graph if and only if it is biregular. See chapter 6 for details.

5. We presented results on the adjacency matrix of $G$-graphs and we calculated the spectrum and energy of some $G$-graphs. See chapter 7.

6. We gave precise examples of infinite $G$-graphs which are locally finite. We presented a result that tells us when an infinite $G$-graph is locally finite. See chapter 7.

## 8.2 Open Questions and Further Works

Generally, much is still to be discovered about $G$-graphs. For instance, in chapter 7 of this thesis, we presented some examples of the spectrum for various $G$-graphs but we were unable to discern any particular behaviour. We also considered infinite $G$-graphs.

One can further investigate the spectrum of a $G$-graph. In fact it is known that, the spectrum of a Cayley graph is directly related to the irreducible characters of the corresponding group [4, 14]. Can we determine the spectrum of a $G$-graph from a property of the group (or/and the generating set)? This will play an important role in its application. Moreover with this, we can easily find the energy of the graph. The energy of a graph is the sum of the absolute values of its eigenvalues. There is interest in which graphs are hypoenergetic and we found two examples.

Another important aspect of a $G$-graph that one can investigate is its *diameter*. The *diameter* of a graph is the longest distance between any two vertices in a graph. This has applications in communication networks[15]. Can we point precisely to the diameter of a $G$-graph from a property of the group, or at least provide a bound? We have not investigated this though.